\documentclass[11pt,notitlepage,twoside,a4paper]{amsart}
 \usepackage{amsfonts}

\usepackage{amsmath,amssymb,enumerate}

\usepackage{epsfig,fancyhdr,color}%,showkeys,amsmidx

\usepackage{amssymb}
\usepackage{amsmath,amsthm}
\usepackage{latexsym}
\usepackage{amscd}
\usepackage{psfrag}
\usepackage{graphicx}
\usepackage[latin1]{inputenc}
\usepackage[all]{xy}
\usepackage[mathcal]{eucal}

\definecolor{NoteColor}{rgb}{1,0,0}

\newcommand{\degree}{\ensuremath{^\circ}}

\renewcommand{\textsc}{\textcolor{red}}

\newtheorem*{theorem 1}{\rm\bf Proposition 1}
\newtheorem*{theorem 2}{\rm\bf Proposition 2}

\theoremstyle{definition}

\theoremstyle{remark}

\def\interieur#1{\mathord{\mathop{\kern 0pt #1}\limits^\circ}}

\title[Euler on spherical geometry]{On the works of Euler and his followers on spherical geometry}

\author{Athanase Papadopoulos}
\address{Athanase Papadopoulos,  Institut de Recherche Math\'ematique Avanc\'ee,
Universit\'e de Strasbourg and CNRS,
7 rue Ren\'e Descartes,
 67084 Strasbourg Cedex, France.} 
 \email{athanase.papadopoulos@math.unistra.fr}

\date{\today}

\begin{document}

\begin{abstract}   We review and comment on some works of Euler and his followers on spherical geometry. We start by presenting some memoirs of Euler on spherical trigonometry. We comment on Euler's use of the methods of the calculus of variations in spherical trigonometry. We then survey a series of  geometrical resuls, where the stress is on the analogy between the results in spherical geometry and the corresponding results in Euclidean geometry. We elaborate on two such results. The first one, known as Lexell's Theorem (Lexell was a student of Euler), concerns the locus of the vertices of a spherical triangle with a fixed area and a given base. This is the spherical counterpart of a result in Euclid's \emph{Elements}, but it is much more difficult to prove than its Euclidean analogue. The second result, due to Euler, is the spherical analogue of a generalization of a theorem of Pappus (Proposition  117 of  Book VII of the \emph{Collection}) on the construction of a triangle inscribed in a circle whose sides are contained in three lines that pass through three given points. Both results have many ramifications, involving several mathematicians, and we mention some of these developments. We also comment  on three papers of Euler on projections of the sphere on the Euclidean plane that are related with the art of drawing geographical maps. 
\bigskip 


\bigskip 

\bigskip 

\noindent  \emph{Acknowledgements.---} The author wishes to thank Norbert A'Campo who taught him several aspects of spherical geometry, the organizers of the International Seminar on History of Mathematics held in Delhi on November  27-28, 2013, and Shrikrishna Dani who gave him the motivation for writing down this article. The work is partially supported by the French ANR project FINSLER. It was finalized during a visit at Galatasaray University (Istanbul) sponsored by a Tubitak 2221 grant. I would also like to thank the two anonymous referees of this article for their careful reading and useful remarks.

\bigskip 
 
\noindent The paper will appear in \emph{Ga\d{n}ita Bh\=ar\=at\=\i } (Indian Mathematics), the Bulletin of the Indian Society for History of Mathematics.

 \end{abstract}
 
\maketitle

\section{Introduction}
The goal of this paper is to bring together some results of Euler and his followers on spherical geometry. By the word ``followers", we mean the mathematicians who benefited from Euler's teaching; some of them were his students, and others were his assistants or young collaborators. Most of them became eventually his colleagues at the Academy of Sciences of Saint Petersburg. 
These works of Euler and his followers contain a wealth of ideas that have not got the attention they deserve from the working geometers. 

The results that we survey can be classified into three categories.

The first set of results concern spherical trigonometry. Euler wrote several papers on that subject, in which he derived a complete set of trigonometric formulae for the sphere. An important contribution in one of the memoirs that we review here is the introduction of the newly discovered methods of the calculus of variations. This allowed Euler to give intrinsic proofs of the spherical trigonometric formulae that are based on the differential geometry of the sphere, unlike the classical proofs where the spherical formulae are derived from the Euclidean,  based on the fact that the lines\footnote{We use the word ``line" in the sense of the prime elements of a geometry, applied to the sphere. This word does not refer to the Euclidean straight lines.} are the intersections of the sphere with the Euclidean planes passing by the origin.

The second category of results that we survey consist of spherical analogues of  Euclidean theorems and constructions.  The idea of examining the analogies between the Euclidean and spherical geometry is very classical, and it can be traced back to the works of the Greek mathematicians Theodosius (second century B.C.), Menelaus (first century A.D.), and Ptolemy (second century A.D.). Understanding and proving the spherical analogues of Euclidean theorems is sometimes not a trivial task, and some of the results obtained by these authors are difficult to prove. We shall mention several examples of such analogies in the works of Euler and his students, and we shall present in detail two of them, together with some developments. The first result is a spherical analogue of a result in Euclid's \emph{Elements} (Propositions 37 and its converse, Proposition 39, of Book I). It characterizes the locus of the vertices of triangles with fixed base and fixed area. The second result is a spherical analogue of a construction by Pappus (Proposition 117 of  Book VII of his \emph{Collection}) of a triangle circumscribed in a circle such that the three lines containing the edges pass through three given points. (In Pappus' Euclidean setting, the three points are aligned.) 

The third category of results that we survey concern maps from the sphere into  the plane. Even if the motivation behind this research is the practical question of drawing geographic maps, the developments are purely mathematical. Euler's main concern in this field is the characterization of maps  from the sphere into the Euclidean plane that preserve specific properties (perpendicularity between the meridians and the parallels to the equator, preservation of area, infinitesimal similarity of figures, etc.) The classical stereographic projections are only one example of such maps.

One important feature of most of the results established by Euler on spherical geometry  is the absence of use of solid geometry in the proofs, and the use of intrinsic methods of the sphere, including polarity theory. This tradition goes back to the work of Menelaus, but this work was not known to Euler.\footnote{Euler and his collaborators were  familiar with the work of Theodosius, but not with that of Menelaus, nor with the later works of the Arabic commentators. Lexell writes in the introduction of  his paper \cite{Lexell-Solutio}: ``From that time in which the Elements of Spherical Geometry of Theodosius had been put on the record, hardly any other questions are found, treated by the geometers, about further perfection of the theory of figures drawn on spherical surfaces, usually treated in the Elements of Spherical Trigonometry  and aimed to be used in the solution of spherical triangles." The work of Menelaus, which in our opinion by far surpasses the one of Theodosius, remained rather unknown until recently. No Greek manuscript survives, but fortunately some Arabic translations reached us.  There exists a German translation of this work, from the Arabic manuscript of Ibn `Ir\=aq \cite{Krause} and there is a  forthcoming English translation from the Arabic manuscript of al-Haraw\=\i \ \cite{RP}; see also \cite{RP1} and \cite{RP2}.}

We have included in this paper some biographical notes on Euler's collaborators, but not on him. There are several very good biographies of Euler, and we refer the reader to the books of Fellmann \cite{Fellmann2}, Spiess  \cite{Spiess1929} and Du Pasquier, 
\cite{Pasquier1927-b}, as well as to the moving tribute (\'Eloge) by Fuss \cite{Fuss-Eloge}, and to the one by Condorcet \cite{Condorcet-Eloge}. The book by Fellmann reproduces a short autobiography which Euler dictated to his oldest son Johann Albrecht.\footnote{Johann Albrecht Euler (1734-1800) was an astronomer and a mathematician. See also Footnote \ref{f:JA}.} We only mention that spherical geometry is one of many fields in which the contribution of Euler is of major importance. The volume \cite{HP} contains articles on several aspects of the works of Euler on mathematics, physics and music theory.

\section{A brief review of the work of Euler on spherical geometry}\label{s:quick}
Before reporting on the works of Euler, let us make some brief remarks on the history of spherical geometry.

 Spherical geometry, as the study of the figures made by intersections of planes with the sphere, was developed by Theodosius. Menelaus inaugurated a geometrically intrinsic study of spherical triangles which is not based on the ambient Euclidean solid geometry, but his work was (and is still) very poorly known, except for some quotes in the work of Ptolemy.  Chasles, in his  \emph{Aper\c cu historique} \cite{Chasles-Apercu} (1837), after mentioning the early works on spherical geometry by Theodosius, Menelaus and Ptolemy, adds the following (p. 236): ``This doctrine [of spherical lines and spherical triangles], which is almost similar to that of straight lines, is not all of spherical geometry. There are so many figures, starting from the most simple one, the circle, that we can consider on this curved surface, like the figures described in the plane. But it is only about forty years ago that such an extension has been introduced in the geometry of the sphere. This is due to the geometers of the North." He then mentions Lexell and Fuss, the two mathematicians who worked at the Saint Petersburg  Academy of Sciences and who had been students of Euler. For instance, Fuss, in (\emph{Nova acta} vol. II and III), studied \emph{spherical ellipses}, that is, loci of points on the sphere whose the sum of distances to two fixed points (called the foci) is constant. Fuss showed that this curve is obtained as the intersection of the sphere with a second degree cone whose centre is at the centre of the sphere. He also proved that if the sum of the two lengths is equal to half of the length of a great circle, then the curve on the sphere is a great circle, independently of the distance of the foci. Some of the works of Lexell and Fuss on spherical geometry  were pursued by Schubert, another follower of Euler, who also worked on loci of vertices of triangles satisfying certain properties. We shall talk in some detail about the works of these mathematicians, and in the last section of this paper, we shall give short biographical information on each of them.

In all the discussion that follows, the distance on the sphere is the angular distance. In other words, we take the model of the unit sphere in 3-space, where the distance between two points is equal to the angle made by the two rays starting at the origin and passing through these points.  A spherical triangle is defined by three points on the sphere together with three lines  joining them, called the sides.\footnote{In working with spherical triangles, it is usually assumed that the length of each of the three edges is smaller than half of a great circle. The reason is that in this case, the edges are shortest arcs between the vertices, and they are the unique such arcs. This hypothesis was already made by Menelaus in his \emph{Spherics}, whose main object is the study of properties of spherical triangles.} The sides of a spherical triangle, being shortest lines, are segments of great circles, that is, they are arcs of intersections of the sphere with planes passing through the center. The angle made by two lines (or great circles) is the dihedral angle made by the corresponding planes.\footnote{We shall stick to this description of spherical geometry, as the geometry of a sphere embedded in 3-space. There is an axiomatic approach of spherical geometry, but this came later on (in the work of Lobachevsky, and then in that of Hilbert).}  The angle at a vertex of a triangle is the angle made by the two lines that contain the two sides that abut at that vertex.

In this section, we give a brief account on the work of Euler on spherical geometry. We shall give below a list of his main memoirs on the subject. The titles should give the reader some idea of their content.\footnote{The reader will notice that the original titles are either in Latin or in French. These were the two official languages at the two academies where Euler worked, the Academy of Sciences of Saint Petersburg and the Academy of Sciences of Berlin. About 80 per cent of Euler's papers are written in Latin and the rest in French. Some of his books are written in German.} Two dates are indicated after each paper; these are  the year  where the paper was presented at the Academy of Sciences (of Saint Petersburg or of Berlin) and the year where the paper was published. The date at which the paper was written is not necessarily the date at which it was presented.\footnote{Some papers were presented several years after Euler's death. A record of these dates is kept on the Euler Archive web site.} For this reason, it is sometimes difficult to make a chronology of the discoveries of Euler, if we want to compare them to those of his students and collaborators on related matters. We know for instance that the memoir \cite{Euler-Variae-T} of Euler, published in en 1797, was presented at the Academy of Saint Petersburg in 1778, that is, 19 years before its publication. The work of Lexell (one of Euler's students) on the same subject was presented in  1781 and published in 1784. The memoirs \emph{Solutio facilis problematis, quo quaeritur sphaera, quae datas quatuor sphaeras utcunque dispositas contingat} and  \emph{Geometrica et sphaerica quaedam} were published in 1810 and 1815 respectively, that is, 27 and 32 years after Euler's death (in 1783). The main reason for this delay is that the Academies of Sciences of Saint Petersburg  and of Berlin received too many of Euler's memoirs and, since the number of pages of each volume the journal, was limited, the backlog became gradually substantial. Starting from the year 1729, and until 50 years after Euler's death, his works filled about half of the scientific part of the \emph{Actes} of this Academy. Likewise, between the years 1746 and 1771, about half of the scientific part of the \emph{M\'emoires}  of the Academy of Berlin consisted in papers by Euler. The memoirs of Euler continued to appear regularly in the \emph{Actes} of the Saint Petersburg Academy, even during the years where he was in Berlin (1741-1766).

 We also know that Euler, in some cases, purposely delayed  the publication of some of his memoirs, in order to leave the primacy of the discoveries to others. This occurred especially in the case where his work and his colleague's work on the same subject were done independently and without knowledge of the others' work. It also happened that sometimes Euler considered that the approach of his colleague was better than than his. A famous example is the attitude he had towards Lagrange's work on the calculus of variations, and we shall elaborate on this in \S \ref{calculus} below. This applies also to the work of some of his his students. Such excessive generosity, which was typical of Euler, is very rare in scientific circles. One can quote here Condorcet (\cite{Condorcet-Eloge} p. 307):
\begin{quote}\small
He never claimed any of his discoveries; and if someone  claimed  something from his works, he was quick to repair such an unintentional unfairness, without too much examining whether  the rigorous equity  demanded from him an absolute withdrawal.
\end{quote}

 Let us now give a list of Euler's main memoirs on spherical geometry.

\begin{enumerate}
\item \label{E1} Principes de la trigonom\'etrie sph\'erique tir\'es de la m\'ethode des plus grands et des plus petits (Principles of spherical trigonometry extracted from the method of maxima and minima),\footnote{We translate by ``method of maxima and minima" the French expression ``m\'ethode des plus grands et des plus petits" This is the name that Euler first used for the calculus of variations.} 1753; 1755 \cite{Euler-Principes-T}.

\item \label{E2}   \'El\'ements de la trigonom\'etrie sph\'ero\"\i dique tir\'es de la m\'ethode des plus grands et des plus petits (Elements of spheroidal trigonometry,  from the method of maxima and minima), 1753; 1755 \cite{Euler-Elements-T}.

\item \label{E3}  De curva rectificabili in superficie sphaerica (On rectifiable curves on spherical surfaces), 1770; 1771 \cite{Euler-Curva1771}.
 
 \item  \label{E4} De mensura angulorum solidorum  (On the measure of solid angles), 1778; 1781  \cite{Euler-Mensura-T}.

\item \label{E5} De repraesentatione superficiei sphaericae super plano (On the mapping of spherical surfaces onto the plane), 1777; 1778 \cite{Euler-rep-1777}.
 
  \item \label{E6} De proiectione geographica superficiei sphaericae  (On a geographic projection of the surface of the sphere), 1777; 1778 \cite{Euler-pro-1777}.
  
   \item \label{E7} De proiectione geographica Deslisliana in mappa generali imperii russici usitata (On Delisle's geographic projection and its use in the general map of the Russian empire), 1777; 1778  \cite{Euler-pro-Desli-1777}.

\item  \label{E8} Trigonometria sphaerica universa, ex primis principiis breviter et dilucide derivata (Universal spherical trigonometry deduced in a concise and clear manner from first principles), 1779; 1782 \cite{Euler-Trigonometria-T}.

\item  \label{E9} Problematis cuiusdam Pappi Alexandrini constructio (On a problem posed by Pappus of Alexandria), 1780; 1783 \cite{Euler-Pappi-T}.

 \item \label{E10} De lineis rectificabilibus in superficie sphaeroidica quacunque geometrice ducendis (On rectifiable lines on spheroidal surfaces drawn geometrically by whatever way), 1785; 1788 \cite{Euler-Lineis}.

\item  \label{E11}   Variae speculationes super area triangulorum sphaericorum (Various speculations on the areas of spherical triangles), 1792; 1797 \cite{Euler-Variae-T}.

\item \label{E13}  Geometrica et sphaerica quaedam (Concerning geometry and spheres), 1812; 1815 \cite{Euler-Geometrica-T}.

    \end{enumerate}

   We did not include in this list any of the memoirs by Euler (or by his students and followers) on astronomy. They also contain results on spherical geometry, but the stress is on a different matter. We are interested here in spherical geometry as a mathematical field which is independent from the science of astronomy. We also did not mention papers by Euler on geometrical problems in 3-space involving spheres (there are several such papers), like the generalizations to spheres of  the famous Apollonian problem of constructing circles that are tangent to three given circles. Such problems, even though they involve spheres, belong to Euclidean geometry and not to spherical geometry.

   Let us make a few comments on the papers in this list.

   The three papers (\ref{E1}) (\ref{E2}) (\ref{E8}) are on spherical trigonometry and we shall elaborate on them in \S \ref{sph}.

      In the paper (\ref{E4}), Euler gives a proof of the formula attributed to Girard for the area of triangles as the angular excess.\footnote{We recall that the area of a triangle with angles $A,B,C$ can be defined as the \emph{angular excess} of this triangle, that is,  
$A+B+C-\pi$. It is well-known that any area function on subsets of the sphere which satisfies the usual additivity properties is, when restricted to triangles, the angular excess (up to a constant factor).} He also gives several other formulae for the area of a triangle in terms of the side lengths. We shall comment on these important formulae in \S \ref{s:area}. In the same paper, Euler uses spherical geometry in order to compute the dihedral angles between the faces of each of the five platonic solids.

   The three papers (\ref{E5}) (\ref{E6}) (\ref{E7}) are on cartography, and more precisely on mappings from the sphere to the Euclidean plane, which produce geographical maps. We shall comment on them in \S \ref{s:maps}.

   The two papers  (\ref{E3}) (\ref{E10}) are on curves on spheres and on spheroids, and we shall report on them in \S \ref{s:alg}.

In the memoir  (\ref{E11}), Euler provides a new formula for the area of a triangle in terms of the lengths of its sides.  We shall elaborate on this question in \S \ref{s:area} below.
In the same memoir, he considers the question of finding the locus of the vertices of 
spherical triangles that have a fixed base and a given area. This is the spherical analogue of a result in Euclid's \emph{Elements}. In the Euclidean case, the required locus is a straight line parallel to the line containing the base of the triangle. See Heath's translation of the \emph{Elements}, \cite{Heath} p. 332-333.

The analogous result in spherical geometry is a theorem which Euler attributes to his student Lexell. We shall report on this important result and its developments in \S \ref{s:Lexell} below.

The memoir (\ref{E13}) contains results on Euclidean plane geometry and their analogues in spherical geometry. One starts with a triangle and three segments, each of them joining one vertex of the triangle to a point on the line containing the opposite side. The question is then to find a condition for the three lines containing the three given segments to intersect in one point.

\begin{figure}[ht!]
\centering
 \psfrag{A}{\small $A$}
 \psfrag{B}{\small $B$}
 \psfrag{C}{\small $C$}
 \psfrag{a}{\small $a$}
 \psfrag{b}{\small $b$}
 \psfrag{c}{\small $c$}
 \psfrag{0}{\small $O$}
\includegraphics[width=.65\linewidth]{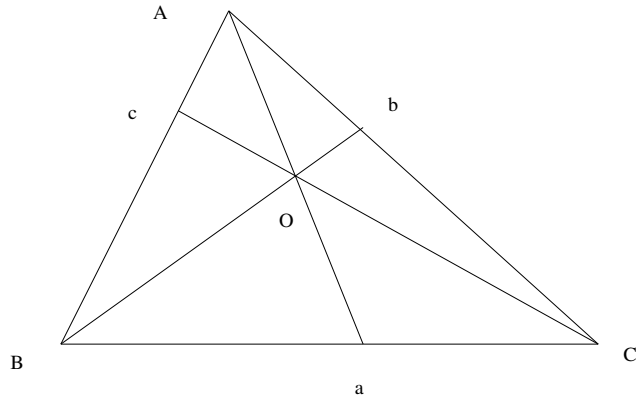}
\caption{\small{ {Equations (\ref{eq:Euler}) and (\ref{eq:Euler2})  give conditions for the lines $Aa, Bb, Cc$ to be concurrent, in the Euclidean and the spherical case respectively.}}}
\label{concourantes}
\end{figure}

 Denoting by $A,B,C$ the vertices of the triangles, $a,b,c$ the three points on the opposite sides and $O$ the intersection point of the three lines (Figure \ref{concourantes}), Euler obtains, in the Euclidean case, the relation  
\begin{equation}
\label{eq:Euler}
\frac{AO}{Oa}\cdot\frac{BO}{Ob}\cdot\frac{CO}{Oc}=\frac{AO}{Oa}+\frac{BO}{Ob}+\frac{CO}{Oc}+2,
\end{equation}
and in the spherical case, 
\begin{equation}
\label{eq:Euler2} \frac{\tan AO}{\tan Oa}\frac{\tan BO}{\tan Ob}\frac{\tan CO}{\tan Oc} =  \frac{\tan AO}{\tan Oa} + \frac{\tan BO}{\tan Ob} + \frac{\tan CO}{\tan Oc} + 2.
\end{equation}
In other words, the condition for the three lines to be concurrent is the same up to replacing the lengths of the sides by the tangents of these lengths. Mathematicians know that such a resemblance is not fortuitous, and it can be explained. We shall see other such analogies below.\footnote{An analogous statement holds in hyperbolic geometry, where the tangent function is replaced by the hyperbolic tangent function. See \cite{PS}, where analogues of results of Menelaus and Euler in spherical geometry are proved in hyperbolic geometry.}

The same memoir contains other results, which are again proved first in the Euclidean case and then in the spherical case. In particular, Euler proves the following converse of the relation (\ref{eq:Euler}): 

Given three segments  $Aa, Bb, Cc$ meeting at the same point  $O$ and satisfying the relation (\ref{eq:Euler}), we can construct a triangle with vertices  $A,B,C$ such that the three points $a,b,c$ are on the opposite sides.

 He then obtains, as  corollaries,  the following relations which are valid in any triangle, in the Euclidean and the spherical cases respectively:

\begin{equation}
\label{eq:Euler22}
\frac{Oa}{Aa} + \frac{Ob}{Bb} + \frac{Oc}{Cc}=1.
\end{equation}

\begin{equation}
\label{eq:Euler221}
\frac{\tan AO}{\tan Oa} + \frac{\tan BO}{\tan Ob} + \frac{\tan CO}{\tan Oc}=1.
\end{equation}

It seems that none of these results, in the Euclidean as well as in the spherical case, was known before Euler. The result (\ref{eq:Euler22}) (in the Euclidean case) is discussed in the paper \cite{Shephard-Euler}.

There is a well-known relation satisfied by three segments  $Aa,Bb,Cc$,  in the situation considered by Euler, so that the three lines they generate are concurrent, namely, \emph{Ceva's Theorem}.\footnote{Giovanni Ceva (1647-1734).}  This is the relation \[\frac{cA}{cB}\cdot \frac{aB}{aC}\cdot\frac{bC}{bA}=1.\]
Euler's relations (\ref{eq:Euler}) and (\ref{eq:Euler22}) are different.

In the memoir (\ref{E9}) \cite{Euler-Pappi-T}, Euler formulates a problem which is more general than a problem which was solved by Pappus, and he solves it in the Euclidean and in the spherical cases. The problem is to construct a triangle which is inscribed in a given circle and such that the sides are contained in  lines passing through three given points. The construction of Pappus concerns the case where the three points are aligned.  We shall mention some developments of this question in \S \ref{s:pappus} below.

\section{Euler's work on spherical trigonometry}\label{sph}

Let us return now to the paper \emph{Principes de la trigonom\'etrie sph\'erique tir\'es de la m\'ethode des plus grands et des plus petits} \cite{Euler-Principes-T}.

As announced in the title, Euler uses in the proofs of the spherical trigonometric formulae the ``method of maxima and minima", which was the name that he used for the calculus of variations.

The application of these methods to spherical geometry is based on the fact that the edge of a triangle is the shortest path between its extremities. The methods of the calculus of variations are then used to find the relations between the six elements of the triangle, namely, the three side lengths and the three angles. In the introduction of the memoir \cite{Euler-Principes-T}, Euler writes:\footnote{The translations from the French are mine.} 
 \begin{quote} \small
Since the sides of a spherical triangle are the smallest lines that one can draw from an angle to another one, the method of maxima and minima can be used to determine the sides of a spherical triangle; and from there we can find after that the relation that remains between the angles and the sides. This is precisely the subject of spherical trigonometry.
\end{quote}

\begin{figure}[ht!]
\centering
 \psfrag{s}{\small $s$}
\psfrag{x}{\small $x$}
\psfrag{y}{\small $y$}
\psfrag{O}{\small $O$}
\psfrag{A}{\small $A$}
\psfrag{M}{\small $M$}
\psfrag{P}{\small $P$}
\includegraphics[width=.85\linewidth]{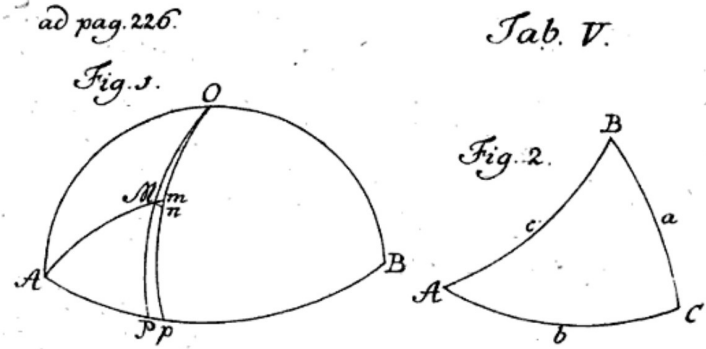}
\caption{\small{The figure is extracted from Euler's \cite{Euler-Principes-T}.}}
\label{Principes1}
\end{figure}

Like in several other texts of Euler, the results in the memoir \cite{Euler-Principes-T} are presented in the form of problems and solutions. 
Problem 1 concerns right triangles on the sphere. Here is the statement:\begin{quote}\small
Given the arc $AP$ on the equator $AB$ and the arc $PM$ on the meridian $OP$, find the shortest line $AM$ which joins the point $A$ to the point $M$ on the surface of the sphere (Figure \ref{Principes1}, left hand side).
\end{quote}

In other words, we want to find the length of a side opposite  to a right angle, knowing the lengths of the two sides which contain this right angle.\footnote{Recall that on the sphere, there exist triangles with more than one right angle.}

In Figure \ref{Principes1}, the arc $AP$ is on the equator, $PM$ on a meridian, and we are searching for the length of the arc $AM$. This arc is produced infinitesimally until the point $m$. Setting $AP=x$, $PM=y$ and $AM=s$, Euler obtains several formulae which solve the problem. Among them, the most well known is probably

\begin{equation}\label{pythagore}
\cos s = \cos x \cos y.
\end{equation}

  This is the so-called \emph{spherical Pythagorean theorem}.\footnote{As the Pythagorean theorem of Euclidean geometry establishes a relation between the three sides of a right triangle, it is also customary in non-Euclidean geometry to call ``Pythagoras theorem" a theorem which establishes a relation between the three sides of a right triangle; cf. \cite{ACP}.} Other formulae are easily deduced from this one, for instance:  

\begin{equation}\label{sinus}
\sin x = \sin  \widehat{MAP} \sin \widehat{AMP}.
\end{equation}

A few remarks may be useful to the reader who might not be aware of the basics of spherical geometry; they will help him make the relation with the well-known formulae of Euclidean geometry.

Let us make the three edges of the triangle $AMP$ very small while keeping the angle $\widehat{AMP}$ constant. We then use the second order Taylor expansions $\cos x\simeq 1-x^2/2$ and $\sin x\simeq x$. The formulae (\ref{pythagore}) and (\ref{sinus}) become respectively  $s^2= x^2+y^2$ and $x=s\sin\widehat{MAP} $, which are well-known formulae of Euclidean geometry:  the first one is the usual Pythagoras Theorem and the second one is just the definition of the sine function. 

These formulae, which are valid for infinitesimal spherical triangles, reflect the fact that infinitesimally, spherical geometry is Euclidean.

Let us also recall that the formulae (\ref{pythagore}) and (\ref{sinus}) assume a certain normalization which comes from
working on a sphere of radius one. For a sphere of a different radius, a certain constant which depends on this radius appears in the formulae.\footnote{In other words, the formula depends on the curvature of the space. In fact, taking this constant into account, one can see that we obtain the Euclidean formulae not only by taking infinitesimal triangles, but also by making the radius of the sphere tend to infinity. This also gives a proof of the well-known fact that a plane is a sphere of infinite radius. The spherical triangles, when the radius of the sphere tends to infinity, become Euclidean.}

\begin{figure}[ht!]
\centering
\includegraphics[width=.85\linewidth]{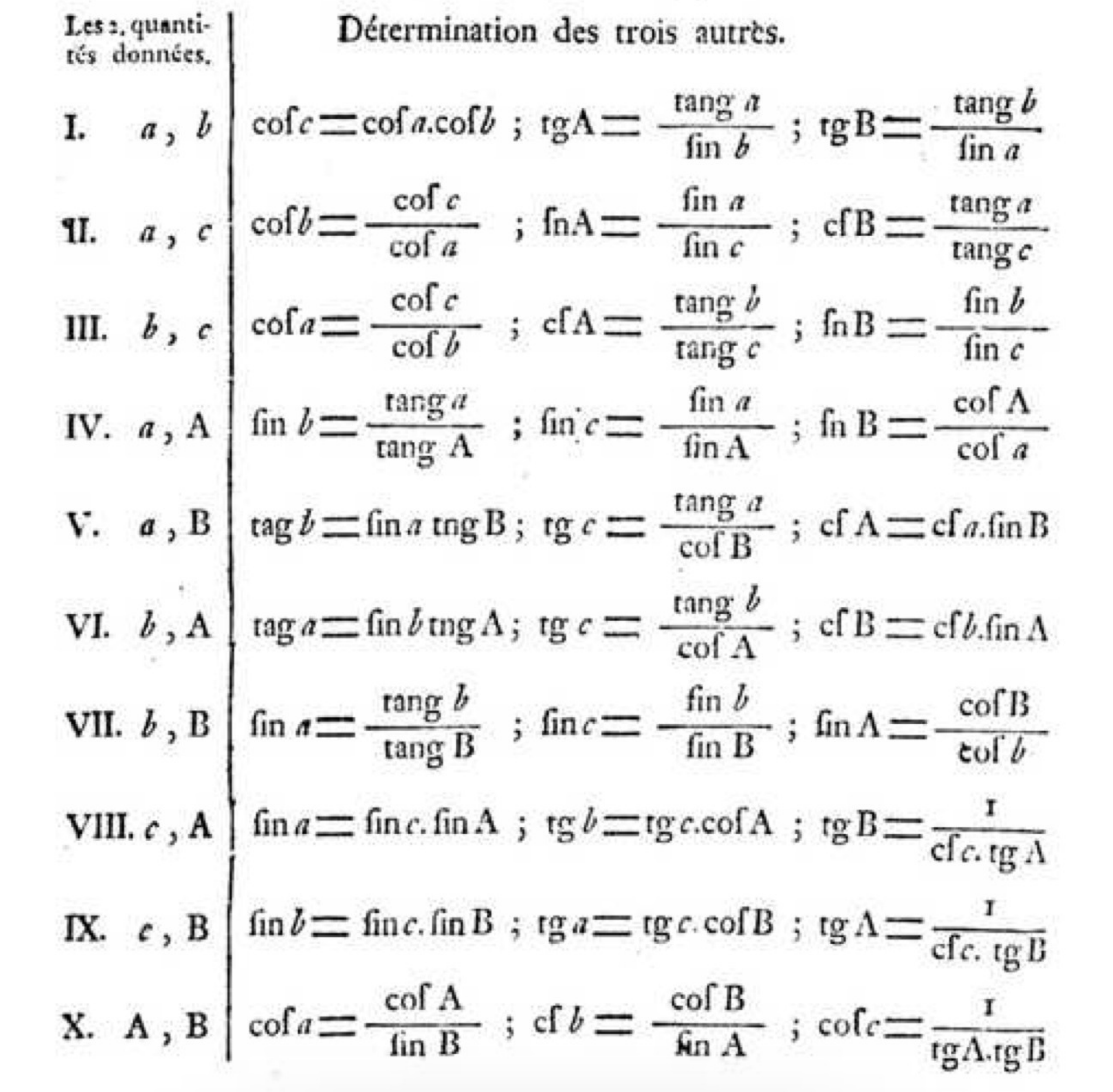}
\caption{\small{A table extracted from Euler's \emph{Principes de la trigonom\'etrie sph\'erique tir\'es de la m\'ethode des plus grands et des plus petits} \cite{Euler-Principes-T}. Euler used to write systematically the formulae with all their transformations, even if this is redundant.}}
\label{table}
\end{figure}

The memoir \cite{Euler-Principes-T} contains 11 fundamental results, with several corollaries for each of them. After the trigonometric  formulae for spherical triangles, Euler establishes  formulae for arbitrary triangles.   
He transforms systematically each formula he obtains in all possible manners, even if the result is redundant; see for instance the table in Figure \ref{table} for a spherical triangle with a right angle at $C$, in which given two quantities in a triangle (for instance, the two sides $a$ ad $b$), one has an expression for a third quantity (the side $c$, or the angle $A$, or the angle $B$). This kind of redundancy is not unusual in Euler's writings and it was probably due to his pedagogical convictions.

In the same memoir, Euler obtains the formula for the area of a spherical triangle as the angular excess (Girard's Theorem which we already mentioned) and  other formulae, like the value of the angles of a triangle in terms of the side lengths, and conversely, values of side lengths in terms of the angles.\footnote{Let us recall that such formulae, giving side lengths in terms of the angles, do not exist in Euclidean geometry, since the angles do not determine side lengths. But there are well-known formulae giving angles in terms of side lengths.} There are also formulae for a side length in terms of the two other side lengths and the angle they contain, and  several formulae of this sort.
  
In using the methods of the calculus of variations in this context, Euler was conscious of the fact that his contemporaries could blame him for resorting to difficult methods to prove results which might be obtained with the classical methods. He writes in 
\cite{Euler-Principes-T}:
\begin{quote}\small
People will no doubt object that it is contrary to the rules of mathematical methods to use the calculus of maxima and minima in order to establish the 
foundations of spherical trigonometry. Besides, that it seems useless to deduce them from other principles, since those which were used up to now are based on elementary geometry whose rigor serves as a rule for all other parts of mathematics. But let me first remark that the method of maxima and minima   acquires here like a new image, since I will show that it may lead us alone to the resolution of spherical triangles. Furthermore, it is always useful to attain by different ways the same truth, since our mind will very likely deduce from it new clarifications.
\end{quote}

Proving the trigonometric formulae of Euclidean geometry using the ambient three-dimensional space is elementary, but it is intellectually non-satisfying because it uses Euclidean geometry whereas the intrinsic geometry of the sphere is non-Euclidean. Euler did not emit such an opinion, but he emphasized the fact that the methods of the calculus of variations that are used to prove the spherical trigonometric formulae have the advantage of being applicable to spaces which are much more general than the sphere or the Euclidean plane. He mentions as examples the spheroidal and conical surfaces.\footnote{Euler calls ``spheroidal surface" a surface obtained by the revolution, along an axis, of an ellipse in three-dimensional Euclidean space. This is also called an ``ellipsoid of revolution", that is, an ellipsoid which has a rotational symmetry about one of its principal axes. We recall that the   Earth is usually approximated by a spheroid.} As in the case of the sphere which he considers in his paper, Euler says that one can analyze the triangles on such surfaces by considering their sides as shortest lines between the vertices. He applies this method in his subsequent memoir \emph{\'El\'ements de la trigonom\'etrie sph\'ero\"\i dique tir\'es de la m\'ethode des plus grands et des plus petits} \cite{Euler-Elements-T}. It is interesting to see  that Euler, in the paper \cite{Euler-Principes-T}, mentions the notion of a triangle \emph{on an arbitrary surface}, as given by three points together with three shortest lines that join them pairwise. Euler writes:
 \begin{quote}\small
 [This research] will provide us with considerable clarifications, not only in spherical trigonometry, but also on the method of maxima and minima. Indeed, since we showed that most of the mechanical and physical problems are solved very promptly by this method, it is most pleasant to see that the same method gives such a great help for the solutions of problems in pure geometry.
  \end{quote}

Let us also quote Euler from  \cite{Euler-Elements-T}:
 \begin{quote}\small 
 Having established the elements of spherical trigonometry on the principle of maxima and minima, my main goal was to determine such a general principle, from which we can deduce the resolution of triangles which are not only formed on a spherical surface but in general on an arbitrary surface.  
 \end{quote}
Without neglecting the practical aspects of these researches, Euler declares in the same memoir:
\begin{quote}\small 
We understand from here that this research may become of great importance, because if the surface of the Earth is not spherical, but spheroidal, then a triangle formed on the surface of the Earth will not be of the kind of which I just talked about.\footnote{The question of the exact form of the Earth was of paramount importance at the time of Euler. We shall come back to this in \S \ref{s:maps} in this paper.} 
\end{quote}
The memoir  \cite{Euler-Elements-T} also contains explicit computations of distances between points on the Earth. 

One of the main ingredients in the memoir \cite{Euler-Principes-T}, in the application of the method of the calculus of variations, is the fact that at the infinitesimal level, the geometry of the sphere is Euclidean. Thus, for instance, in the solution of the problem of finding the shortest line between two points on the sphere (Problem VI, p. 235 of \cite{Euler-Principes-T}), and using the notation of Equations (\ref{pythagore})  and (\ref{sinus}) and denoting the infinitesimal arc arcs $Mm$, $Pp$ and $mn$ by $ds$, $dx$ and $dy$ respectively, the length element $ds$ is written as 
\[ds = \sqrt{dx^2+dy^2\sin^2x}.\]
This formula is obtained by applying the Euclidean Pythagorean formula to an infinitesimal right triangle. The problem is then to search for conditions on the length $\int \sqrt{dx^2+dy^2\sin^2 x}$ of a path joining two points so that this length is minimal.
This leads to a differential equation and the trigonometric formulae are obtained by solving such differential equations.

After finding some basic trigonometric formulae using the calculus of variations, other formulae are proved by algebraic transformations. The formula for the area of a triangle uses an integration.

Beyond the obtention of trigonometric formulae, Euler discovered, with his use of the calculus of variations, a version of the intrinsic differential geometry of surfaces, a theory which Gauss and then Riemann were about to develop more thoroughly, in the century that followed. One may add here that Euler investigated, long before Gauss, the differential geometry of surfaces embedded in 3-space. In 1772, he studied developable surfaces and he obtained a criterion for two surfaces to be applicable on each other,\footnote{In modern terms, these are surfaces that are locally isometric.} a criterion which was rediscovered in 1828 by Gauss. In 1763, Euler started a thorough study of curvature of embedded surfaces. In 1767, he found an expression of the curvature in terms of the product of principal curvatures.

The memoir   \cite{Euler-Trigonometria-T} is another complete treatise on spherical trigonometry. There is however a major difference with the preceding one: Euler does not use any more the variational methods to prove the trigonometric formulae, and instead he works them out in a classical manner, deriving them from solid geometry. In fact, he applies a spherical triangle onto a tangent plane passing through one of the vertices, using radial projection from the origin of the sphere. Right at the beginning of the paper, he notes that the choice of the notation $A,B,C$ and $a,b,c$ for the angles and sides opposite to them leads to a set of symmetric formulae; for instance, we have the following (extracted from his paper):
\[
\includegraphics[width=.90\linewidth]{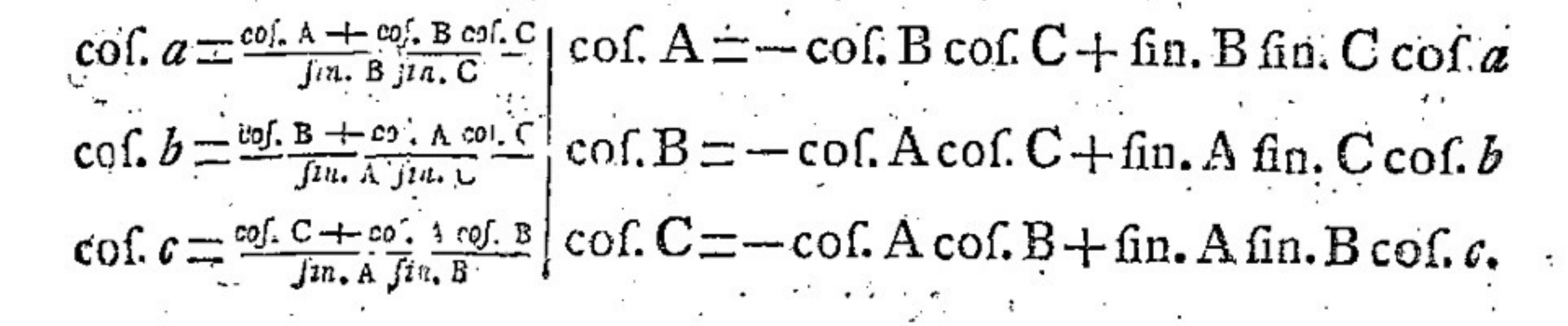}
\]
\[
\includegraphics[width=.90\linewidth]{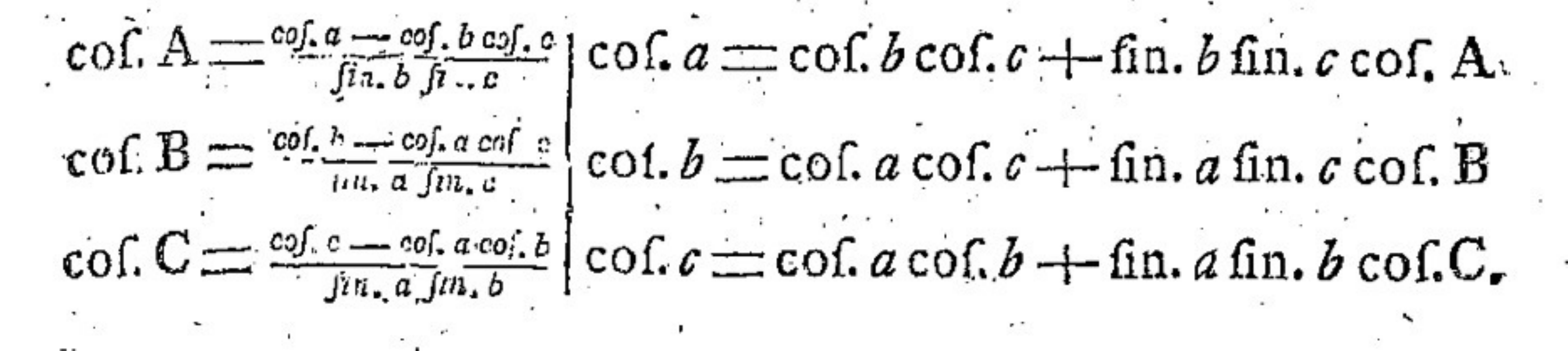}
\]

This is also typical of Euler's writing: to publish different proofs of the same result.

In the paper \cite{Euler-Trigonometria-T}, Euler highlights a rule for obtaining new formulae out of known ones, and this rule is no other than the duality theory in spherical geometry. We know by the way that until recently, the theory was attributed to Vi\`ete\footnote{See e.g. \cite{Chasles-Apercu}.}, but we know now that it was discovered by the Arabic mathematicians of the tenth century, see e.g. the historical work of Marie-Th\'er\`ese Debarnot \cite{Debarnot}. In fact, Euler, in the paper \cite{Euler-Mensura-T} to which we refer in \S \ref{s:area}, uses the notion of polar triangle, which is one way of explaining the symmetry of the formulae. He establishes, using solid geometry, the following three basic formulae,
\[\frac{\sin C}{\sin c}=\frac{\sin A}{\sin a},
\]
\[
\cos A \sin c = \cos a \sin b - \sin a \cos b \cos C,
\]
and
\[
\cos c= \cos a \cos b +\sin a\sin b\cos C,
\]
and he then derives from them the rest of the trigonometric formulae, by performing algebraic transformations. 
   
 Several authors after Euler worked out, like him (and they were probably inspired by him), a complete system of trigonometric formulae based on very few starting formulae; see e.g. Lagrange \cite{Lagrange-Solution-T}  and de Gua \cite{Gua}.

\section{Euler and the calculus of variations}\label{calculus}

In this section, we give, for the convenience of the reader, a short introduction to the calculus of variations, since this theory plays an important role in the work of Euler on spherical trigonometry. In some sense, this field is the study of certain extrema 
  (maxima and minima) of functions defined on some space, usually infinite-dimensional.\footnote{Even though a precise definition of an ``infinite-dimensional space" was not known in the times of Euler.} A typical example of space on which the function is defined is a space of maps. For instance, it could be a space of paths joining two given points on a surface, and we may ask that the paths are differentiable, or piecewise differentiable, etc.\footnote{A path in a space $E$ is regarded here as a map from the interval $[0,1]$ into $E$.} The domain considered may also be a space of surfaces, etc. A function defined on a space of functions is often called a  \emph{functional}. Thus, the object of the calculus of variations is to determine the extrema of functionals. For instance, Morse theory, which was born in the twentieth century and whose object is the study of submanifolds of a differentiable manifold, uses widely the techniques of the calculus of variations, the space of submanifolds being here the underlying infinite-dimensional space. In a certain sense, the calculus of variations is a generalization of differential calculus where one studies extrema of functionals instead of extrema of functions. The idea of generalizing the notion of function to that of functional may seem trivial, but practically, the new theory involves a wealth of new delicate ideas as well as powerful tools.  
  
  The functionals of which we study the extrema in the calculus of variations are usually defined by integrals. A typical example of a space on which the functional is defined is the space of paths between two fixed points and where the functional is the length function. A point where an extremum of this function is attained is a \emph{geodesic}.\footnote{In all this paper, the term \emph{geodesic} means a curve in some space (which is generally a surface) that is shortest among all the curves joining its endpoints.} The determination of geodesics, for instance on a surface, is a fundamental problem in geometry. To give an elementary example, let us recall that in the Euclidean plane, using cartesian coordinates, the length of a curve defined as the graph of a function $f$ defined on the interval $[0,1]$ is given by the integral
  \[\int_0^1\sqrt{1+(f'(x))^2}dx.\]
One can see on this example that some hypotheses must be satisfied by the paths on which this functional is defined, e.g. that it should be differentiable (but there are more restrictive hypotheses).
  
  In fact, one of the first papers that Euler wrote concerns geodesics, \emph{De linea brevissima in superficie quacunque duo quaelibet puncta jungente} (On the shortest line joining two arbitrary points on a surface) \cite{Euler-Brevissima-T} . He wrote it in 1728, at the age of 21. The problem considered in this paper is typical of the calculus of variations. Working in a more general context, Euler gives in this paper a necessary and sufficient condition for a path to be an extremum of the length functional,  a condition which can be expressed by a differential equation which is called today the \emph{Euler-Lagrange equation}.\footnote{Note that Euler also used the Euler-Lagrange equations in his memoir \cite{Euler-Principes-T}, \emph{Principes de la trigonom\'etrie sph\'erique tir\'es de la m\'ethode des plus grands et des plus petits}, which we already mentioned.}

Although the calculus of variations, as a field, was born with Euler, there are problems which belong to that field which were studied prior to Euler. Let us mention a few of them.

  Newton, in Book II of the \emph{Principia}, considered such problems, in the setting of the motion of objects in fluids. He searched for a geometric characterization of a surface of revolution moving at a constant speed in a fluid, in the direction of its axis, such that the air resistance exerted on this surface is minimal.\footnote{\emph{Prinicipia Phil. naturalis mathematica} (1687), sec. II, Prop. 35. Cf. Kline \cite{Kline-T} p. 257 for an exposition.} The british clockmaker Henry Sully (1680-1729), who owned a clock company in Versailles, proposed an isochronal pendulum\footnote{\label{iso}A curve is said to be isochronal if it takes the same time for a body sliding on it without friction and starting at any point to get to the bottom, under the sole force of gravity. The cycloid is an example of such a curve.} which Euler studied later in his article
\cite{Euler-tauto-T}, as a problem of calculus of variations. 
 Johann Bernoulli (1667-1748) and his brother Jakob (1654-1705), as well as other mathematicians before Euler, worked on a problem which is typical of the calculus of variations, namely, the \emph{brachistochrone}, or the ``curves which realizes the minimal time". This is a curve situated in a vertical plane on which a point, submitted to a constant gravitational field of the Earth, sliding without friction, attains the lower point in a minimum amount of time.  This curve is not a line (unless the lower point lies on the same vertical).\footnote{In other words, even though the straight line between two points is length-minimizing, it is not time-minimizing.}
 Galileo, in his   \emph{Discourses and Mathematical Demonstrations concerning Two New  Sciences}  (Theorem XXII, Proposition XXXVI), searched for a characterization of this curve, and he concluded (erroneously) that it must be an arc of circle.\footnote{The reasoning of Galileo, even if its false, is subtle. Let us recall it. There are three steps. The first step consists in showing that given a vertical line $d$  and a point $P$ which is not on $d$ and which is subject to a constant gravitational field, the line segment on which $P$ can slide in order to attain in a minimum time a point on $d$ is the segment joining $P$ to a point $P'$ of $d$ making with this line an angle of $45^{\mathrm{o}}$. The second step consists in proving that if we take a point $P''$ on an arc of circle joining $P$ to $P'$, then if the point $P$ slides over the broken line $PP''-P''P'$,  it attains the point $P'$ more rapidly than if it slides on the first straight line. The third step (the one which is not correct) consist in deducing from this reasoning that the line that we are seeking for is an arc of a circle joining  $P$ to $P'$.} The brachistochrone problem is a typical problem in the calculus of variations. The space on which the function is defined is the space of paths from $P$ to $P'$. The functional that is minimized in this setting is the time of descent, without friction.  
 
 The question of determining the brachistochrone between two points has a long history. It was proposed by Johann Bernoulli in 1696. Bernoulli, who knew the answer, submitted the problem as a challenge to his colleagues.\footnote{Cf. \emph{Acta Eruditorum Lipsiensia}, 1691, p. 269; \emph{Opera} I, p. 161. The problem was formulated by Bernoulli as follows: 
 \begin{quote}
 \emph{Problema novum ad cujus solutionem mathematici invitantur : Datis in plano verticali duobus punctis $A$ et $B$, assignari mobili $M$ viam $AMB$, per quam gravitate sua descendens, et moveri incipiens a puncto $A$, brevissimo tempore perveniat ad alterum punctum $B$} (A new problem which the mathematicians are invited to solve: Given, in a vertical plane, two points $A$ and $B$, determine the trajectory $AMB$ described by a mobile point $M$ which descends by virtue of its weight, starting at $A$ and arriving in the shortest possible time at the other point $B$).
 \end{quote}
  Bernoulli first gave his colleagues six months for the solution, but Leibniz asked him to give more, so that the foreign mathematicians could also contribute to the problem.} The following year, Newton, Leibniz, L'Hospital, Johann Bernoulli himself and his brother Jackob, published a solution, using different methods.\footnote{One can recall here that the atmosphere was very tense between the two Bernoulli brothers, who were constantly in competition with each other. Concerning the problem of the brachistochrone, Johann Bernoulli addressed to the Royal Academy of Sciences a sealed letter containing a solution of the problem, ordering that the letter could not be opened before his brother Jakob makes his solution known. It appeared later on that the solution proposed by Johann Bernoulli was not correct, whereas the one given by his brother Jakob was correct, cf. \cite{Strauch}. It is also considered that Johann Bernoulli tried to plagiarize his brother regarding this question; cf.   \cite{Smith-T} p. 645, which contains an excerpt of the article of  Johann Bernoulli and a small commentary on the history of the question.} Jakob  Bernoulli  showed that this curve was a  cyclo\"\i d.\footnote{A cyclo\"\i d, known also as ``Aristotle's wheel", is the trajectory of a point of a circle moving without friction on a straight line. This curve was well known at the time of Euler, and it was already used by Huygens and others.} The method consisted in varying the curve at a given point and comparing it to the curve thus obtained. It was only after the death of Jakob Bernoulli (1705) that Johann Bernoulli admitted that the methods used by his brother were correct.  In 1719, he published an exposition together with a commentary of these methods,  where he extracted an equation which could be used as the basis of a general theory (cf. \cite{Fraser} p. 110 and  117). In 1734, Euler generalized the problem of the brachistochrone taking into account the resistance of the ambient space,  to questions of extrema of quantities that can be different from time. 
 
 It is fair to note here that although there were problems that were considered before Euler and that are typical of the calculus of variations, the methods used to solve them were geometric and not variational. For instance, geodesics on surfaces were characterized by properties of the osculating planes at these points. They did not involve the analytical methods which were introduced by Euler and Lagrange and which became the methods of the calculus of variations.

The article \emph{Constructio linearum isochronarum in medio quocunque resistente} \cite{Euler-Constructio-T} of Euler,  which concerns the construction of an isochronal curve\footnote{See Footnote \ref{iso}.} in a resistant medium, is his first published work (he was 18 years old); thus, it is by a work on the calculus of variations that Euler started his scientific production. This subject pursued Euler during the rest of his life, and several other questions concerning motion in a resistant medium were treated later in a much more exhaustive manner in his \emph{Mechanica analytica} (1736). The article \cite{Euler-Constructio-T}, as well as the article \cite{Euler-Brevissima-T} which we already mentioned and which concerns geodesics on surfaces,  were also motivated by questions which were asked to Euler by Johann Bernoulli. Euler acknowledges, in the introduction to \cite{Euler-Brevissima-T}, that it was Johann Bernoulli who encouraged him to work on the question of determining the shortest line between two arbitrary points on a convex surface. He writes that ``to solve the problem, it is now necessary to use the \emph{method of maxima and of minima}." He also announces in the introduction that Johann Bernoulli had already found an equation for the shortest line joining two points on an arbitrary surface. In the same article, Euler considers cylinders (whose base is not necessarily circular), cones with circular bases, and surfaces of revolution which are not necessarily convex.

Among the other writings of Euler on the calculus of variations, we can mention his \emph{Problematis isoperimetrici in latissimo sensu accepti solutio generalis} (General solution of the isoperimetric problem considered in the largest sense) \cite{Euler-Problematis-T}, presented to the Academy of Saint Petersburg in 1732, his \emph{De linea celerrimi descensus in medio quocunque resistente} (On the curve of fastest descent in an arbitrary resistant medium) \cite{Euler-Linea-T}, presented in 1734, his \emph{Curvarum maximi minimive proprietate gaudentium inventio nova et facilis} (A new and easy method for finding curves which have a maximal or a minimal property) \cite{Curvarum-T} presented in 1736, his  \emph{Solutio problematis cuiusdam a Celeberrimo Daniele Bernoullio propositi} (Solution of a certain problem proposed by the very famous Daniel Bernoulli) \cite{Euler-Sol-T}, presented in 1738, and the work \emph{Methodus inveniendi lineas curvas maximi minimive proprietate gaudentes, sive solutio problematis isoperimetrici latissimo sensu accepti} (Methods for finding curves  which possess a certain property of maxima and minima, or the solution of isometric problems in the largest sense) \cite{Euler-Methodus-T}, of which Euler wrote a first version in 1740 and which was published in the form of a book in 1744. With all these publications, the subject of the calculus of variations was transformed progressively, by Euler, from the study of particular cases into a general theory.

The work \cite{Euler-Methodus-T} contains a general formulation of the problem of the calculus of variations as well as an exposition of the general methods  to solve it. Euler gives, as an illustration, a list of a hundred particular cases classified in seven categories to which the method applies. The publication of this work is considered as marking the beginning of the calculus of variations as a new branch in mathematics.  We also find there a new and elegant solution of the brachistochrone problem.  This problem is formulated as the one of finding the function $y=y(x)$ which minimizes the quantity \[\int_a^b Z(x,y,y',\ldots,y^{(n)})dx\]
where  $Z$ is a function which, in modern terms, is called the Lagrangian. The work also contains applications of the calculus of variations to physics, and in particular, a discussion of the  principle of least action. Certain problems in geometric optics are also naturally formulated there in terms of the calculus of variations.\footnote{Let us recall that Euler wrote several memoirs on optics.}  

  Carath\'eodory, who was the editor of the volume of Euler's \emph{Opera Omnia} that contains the \emph{Methodus} \cite{Euler-Methodus-T}, describes this work as being  ``one of the most beautiful mathematical works ever written." Later on, Euler wrote his \emph{Analytica explicatio methodi maximorum et minimorum} (Explanation of the analytical method for maxima et les minima), presented to the Academy of Berlin  in 1756 \cite{Euler-Analytica-T} and his \emph{Elementa calculi variationum} (Elements of the calculus of variations)  \cite{Euler-Elementa-T} (written eight days after the first one), using the name that is known today for that theory.\footnote{Euler had already used this name in his correspondence. The expression ``calculus of variations" refers to the introduction and the use of an increment, or ``variation" $\delta$, in the approach to the problem. Both terms, ``method of the maxima and minima" and ``calculus of variations" denote the same field, but they stress different aspects of it. The first expression, despite the fact that it involves the word ``method", designates the goal to attain, which is to find or to prove the existence of a solution defined by an extremal property. The second name, ``calculus of variations", designates the method, viz., taking a small variation.} Since Euler attributes the invention of this theory to Lagrange, let us say a few words on the works of the latter.  

%\begin{center}
% \includegraphics[width=11cm]{./Photos/Lagrange-2.pdf} \\
% \vspace{0.5cm}
%\end{center}
%\noindent
%\begin{minipage}{\textwidth}
%\begin{center}
% \vspace{0.8cm}
% \includegraphics[width=11cm]{./Photos/Lagrange-2.pdf} \\
%% \textit{\small Henri Poincar\'e.}
% \vspace{0.8cm}
%\end{center}
%\end{minipage}
 
In a letter dated August 12, 1755, Lagrange,\footnote{Ludovico de la Grange (1736-1813), better known under the name Joseph-Louis Lagrange, was appointed mathematics professor at the Turin Royal Academy of artillery two months after he wrote that letter to Euler. Youschkevitch and Taton, in the introduction to Volume V of the Series IVA of the \emph{Opera Omnia}, which contains the correspondence between Euler and Lagrange, write (p. 34) that ``it is only in the person of Joseph-Louis Lagrange that Euler found, in 1754, an interlocutor who has an intellectual status which is comparable to his own, and who, at the same time, is a scientist who, since his first steps, turned out to be the most gifted and the most original successor of his work. [...] Their exchanges, marked by a harmony in the intellectual interests, practically encompass all the branches of pure and applied mathematics." We can recall by the way that Lagrange was the successor of Euler at the Academy of Berlin, in 1766,  after the return of the latter to the Academy of Sciences of Saint Petersburg. On the other hand, Euler was very generous with Lagrange, but the latter did not always express gratefulness to Euler. It suffices to skim through the correspondence between Lagrange and d'Alembert in order to see the condescending and sometimes derisive  tone which Lagrange uses when he talks about Euler.} who was only nineteen years old, explained to Euler, who was thirty years older than him, new ideas he discovered on the calculus of variations. These ideas simplified the works of Euler in the sense that they avoided delicate geometric arguments and replaced them by analytic arguments which lead directly to what was called later on the Euler-Lagrange equation. Lagrange had already written to Euler on June 28, 1754, on another subject,\footnote{Lagrange talked in this letter about the analogy between the powers of the binomial $(a+b)^m$ and the differentials of the product $d^m(xy)$, and he asked Euler about his opinion on this analogy whose discovery he considered as his own; cf.  \emph{Opera Omnia}, S\'erie IVA vol. V, p. 361 and the commentary by Youschkevitch and Taton on p. 35  of the same volume.} but this first letter remained without response. Euler responded to the second letter on September 6 of the same year, expressing to Lagrange his great joy in learning these discoveries and congratulating him of ``carrying his [that is, Euler's] theory of maxima and minima to its highest degree of perfection". The following year, Euler presented to the Academy of Sciences of Berlin the new methods of Lagrange. Lagrange published his results in 1760 in an article titled \emph{Essai d'une nouvelle m\'ethode pour d\'eterminer les maxima et les minima des formules int\'egrales ind\'efinies} (An essay for a new method for determining the maxima and the minima of indefinite integral formulae), which appeared in the \emph{Miscellanea taurinensia}, a journal from Turin. The article contains solutions to problems in the calculus of variations which are more general than those that Euler studied. For instance, Lagrange considers the brachistochrone problem with variable extremities.\footnote{In his letter to Euler dated November 20, 1755, Lagrange had already presented the solution of the brachistochrone problem in which the origin is fixed and the extremity varies on a given curve.} It is in this article that the so-called Lagrange multipliers appear for the first time.

 We know that Euler waited several years before publishing some of his papers on the calculus of variations, because he wanted the primacy of the discovery to go to Lagrange. This is apparent in the correspondence between the two men (which lasted 21 years); cf. \cite{Euler-correspo1-T}. We learn in particular (p. 421ff.) that at the birth of this theory, the works of Euler were shading off those of Lagrange, despite the efforts of Euler to highlight those of his young colleague. The reason was mainly that the style of Lagrange  was too concise,\footnote{Euler, who was prolific, was an admirer of the concise style of Lagrange. Youschkevitch and Taton,  in the introduction to Volume V of the Series IVA of the \emph{Opera Omnia},  write (p. 35): 
``Euler used to admire not only the profoundness of the researches of Lagrange but also their surprising concision."} which made his writings difficult to read. In a letter addressed to Euler on October 28 (?) 1762,\footnote{\emph{Opera Omnia}, S\'erie IVA vol. V, p. 446.} Lagrange summarizes the situation in the following terms:
\begin{quote}\small  
Having learned, by one of your letters of 1759, that you rate my method of  \emph{de maximis et minimis} enough highly so as to extend it and improve it in a special treatise, I thought I would had to delete entirely the one which I had already finished on the subject and I limited myself to the simple exposition of the principles in a memoir which I tried to make as concise as possible. I even decided to write this memoir only because you gave me in that letter the honor of informing me that you will not publish your work before I publish mine. I look forward to taking advantage of the new lights that you have spread on such a difficult matter. 
\end{quote}

The same year, Euler exposited the methods of  Lagrange, using the notation of the latter. In 1788 (five years after the death of Euler) Lagrange's \emph{M\'ecanique analytique} appeared. This is probably his most important piece of work. It contains the foundations of the new mechanics where Euler's geometrical methods are completely replaced by analytic methods.

 \section{On the locus of vertices of triangles having the same base and the same area}\label{s:Lexell}
 
   We now return to spherical geometry. We consider, in this section and the next one, works on  construction problems of triangles. We start with a work of Lexell. Section \ref{s:students} of the present paper contains some biographical elements on Lexell.

Lexell published a memoir titled \emph{Solutio problematis geometrici ex doctrina sphaericorum} (The solution of a geometrical problem according to the spherical doctrine) \cite{Lexell-Solutio}, in which he studied the locus of the vertices of triangles with a given base and fixed area. We mentioned this problem in the introduction.
 
\begin{figure}[ht!]
\centering
 \psfrag{A}{\small $A$}
 \psfrag{B}{\small $B$}
 \psfrag{C}{\small $C$}
 \psfrag{O}{\small $O$}
 \psfrag{M}{\small $M$}
 \psfrag{P}{\small $P$}
 \psfrag{Q}{\small $Q$}
 \psfrag{Z}{\small $Z$}
 \psfrag{V}{\small $V$}
\includegraphics[width=.70\linewidth]{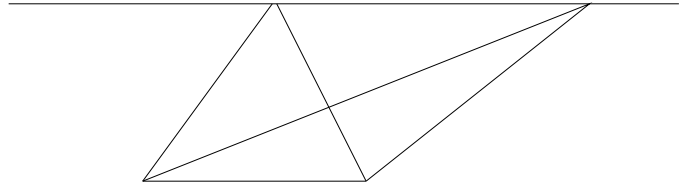}
\caption{\small {In the Euclidean plane, the locus of the triangles having a fixed area and the same base is a line parallel to the base.}}
\label{Eucl}
\end{figure}

We recall that Euclid in the \emph{Elements} solves the analogous Euclidean problem, and the solution there, that is, the locus of the vertices of a triangle of fixed area and fixed base, is a line which is parallel to the base (Figure \ref{Eucl}).

 Euclid calls two triangles having the same area \emph{equal}. The reason is that such triangles can be dissected in such a way that the pieces can be reassembled to lead to the same triangle. The same definition of ``having the same area" applies to spherical triangles.   We recall by the way that Euclid did not use the absolute notion of \emph{area}, which was beyond his means.\footnote{We know that in order to define properly the notion of area in Euclidean geometry, one has to develop measure theory, which is a complicated theory: one quickly encounters there the notion of non-measurable set, and so on.}

 In the spherical case, the solution of the problem is not a line but a small circle which is \emph{not} equidistant to the base of the triangle.\footnote{We recall that on the sphere, a \emph{circle} is a locus of points that are at a certain distance (the \emph{radius} of the circle) from a give point (the \emph{center} of the circle). A circle of maximal radius is a \emph{great} circle. The other circles are called \emph{small} circles. The great circles are the ``lines" of a sphere, that is, its geodesics. (We already mentioned that in the study of the geometry of the sphere, the sides of spherical triangles are segments of great circles.)  In Euclidean geometry, the set of points that are equidistant from a given line is a union of two lines. In spherical geometry, this set is a union of two small citrcles (provided the distance is small enough -- smaller than the radius of the sphere. A small circle is also a component of an equidistant locus from a great circle.}

Let us report on the solution of Lexell. We use his notation.

\begin{figure}[ht!]
\centering
 \psfrag{A}{\small $A$}
 \psfrag{B}{\small $B$}
 \psfrag{C}{\small $C$}
 \psfrag{O}{\small $O$}
 \psfrag{M}{\small $M$}
 \psfrag{P}{\small $P$}
 \psfrag{Q}{\small $Q$}
 \psfrag{Z}{\small $Z$}
 \psfrag{V}{\small $V$}
\includegraphics[width=.80\linewidth]{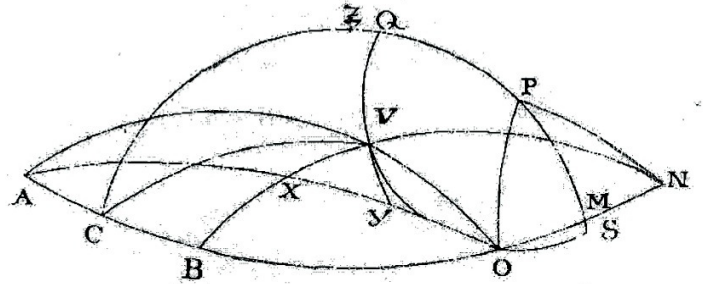}
\caption{\small {The construction of the spherical locus that corresponds to the Euclidean locus of Figure \ref{Eucl}. (Extracted from Lexell's paper \cite{Lexell-Solutio}).}}
\label{Lexell1}
\end{figure}

We are looking for the locus of the vertex $V$ of a triangle $AVB$ with  fixed base $AB$ and with fixed area 2$\delta$. (Thus, the angle sum of the triangle $AVB$ is $180\degree + 2\delta$.) The locus that we are seeking is the segment of small circle $QVO$ represented in Figure \ref{Lexell1}. 

The construction is as follows.

Let $O$ be the second intersection point of the lines $AB$ and $AV$. (In other words, $O$ is the antipode of the point $A$.) Let $C$ be the midpoint of $AB$ and let $CZM$ be a line (great circle) perpendicular to $AB$ at $C$, where $M$ is the second intersection point of this line with the line $AB$. (Thus, $M$ is the antipode of $C$.) From the point $O$, consider the line $OP$, where $P$ is on $CZS$ and where the angle $\widehat{POM}$ is equal to $90\degree -\delta$. Finally, taking $P$ as a pole we draw the small circle  at distance $PO$; this circle passes through $V$ and  it intersects the line $CZM$ at a point $Q$. The arc $QVO$ is the locus that we are looking for. Let us note that the reasoning shows that the locus is a small circle that passes by the antipode of $A$, and therefore the antipode of $B$ as well.

During the construction, Lexell makes the following remarks (\cite{Lexell-Solutio}, \S 6 and 10):

--- For any vertex $V$ of a triangle satisfying the requirements, the arc of great circle $AV$ intersects the great circle containing the side $AB$ at the same point $O$, and we have $MO=BC$,  which is half of the distance $AB$.

--- There is a limiting case where the point $V$ is at $O$. In this case, the triangle $AVO$ becomes equal to a bigon $AXOBA$. 

--- The segment $PO$ is perpendicular to the line $AXO$. 

--- The locus $OVQ$ is always a segment of a \emph{small} circle, except in the case where  $CB$ is a quarter of a circle, and in this case the locus $OVQ$ is a segment of great circle.

It is interesting to compare the solutions of Lexell and of Euler.\footnote{Euler's memoir was published posthumously in 1797, 16 years after the one of Lexell. The two memoirs, of Lexell and of Euler, where written in the same year.} In the memoir \cite{Euler-Variae-T}, Euler declares that the impulse for studying this problem came to him from the result of Lexell.

In \S 16 of \cite{Euler-Variae-T}, after giving the formula for the area of a spherical triangle in terms of its side lengths, Euler writes (translation by J. C.-E. Stern):
\begin{quote}\small
The occasion for me to start pondering on this came from a theorem concerning all the spherical triangles having the same area raised upon the same base, brought into light by the famous Professor Lexell, who acutely demonstrated that all the vertices of these triangles always lie on some small circle of the sphere, which most elegant property may be derived not without many detours from our theorem.
\end{quote}

Euler's proof of Lexell's theorem starts with a few lemmas on quadrilaterals on the sphere, which he calls \emph{parallelograms}.

Consider on the sphere two small circles which are situated on different sides of a given great circle and at the same distance. (Note that the two small circles have the same poles and therefore they have equal lengths.)

\begin{figure}[ht!]
\centering
 \psfrag{M}{\small $M$}
 \psfrag{E}{\small $E$}
 \psfrag{N}{\small $N$}
 \psfrag{m}{\small $m$}
 \psfrag{n}{\small $n$}
 \psfrag{e}{\small $e$}
 \psfrag{f}{\small $f$}
 \psfrag{O}{\small $O$}
 \psfrag{A}{\small $A$}
  \psfrag{B}{\small $B$}
 \includegraphics[width=.55\linewidth]{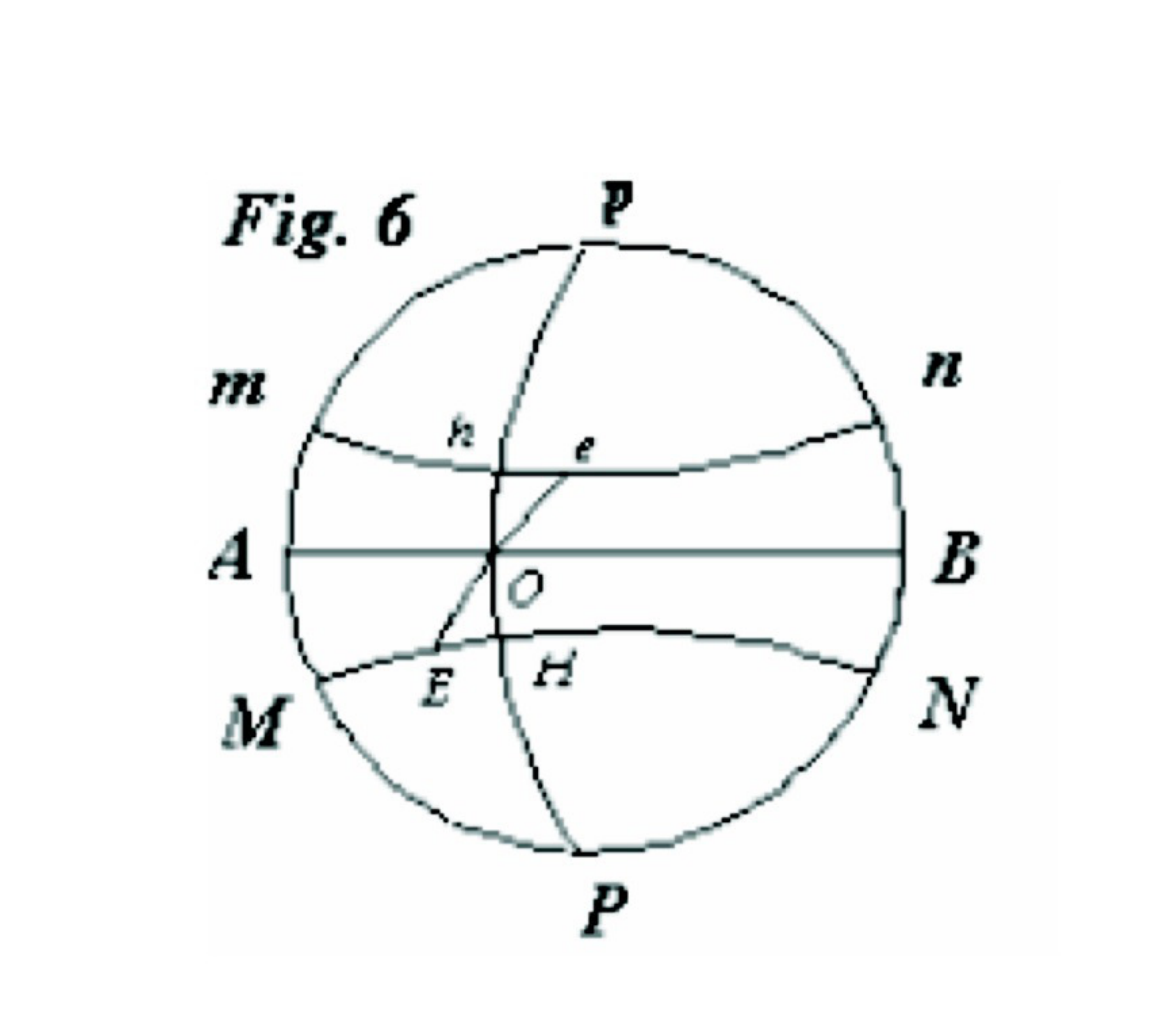}
\caption{\small {Extracted from Euler's paper \cite{Euler-Variae-T}.}}
\label{para2}
\end{figure}
\begin{figure}[ht!]
\centering
 \psfrag{z}{\small $\zeta$}
 \psfrag{x}{\small $\xi$}
 \psfrag{e}{\small $e$}
 \psfrag{f}{\small $f$}
 \psfrag{E}{\small $E$}
 \psfrag{F}{\small $F$}
  \psfrag{m}{\small $m$}
 \psfrag{n}{\small $n$}
  \psfrag{M}{\small $M$}
 \psfrag{N}{\small $N$}
 \includegraphics[width=.55\linewidth]{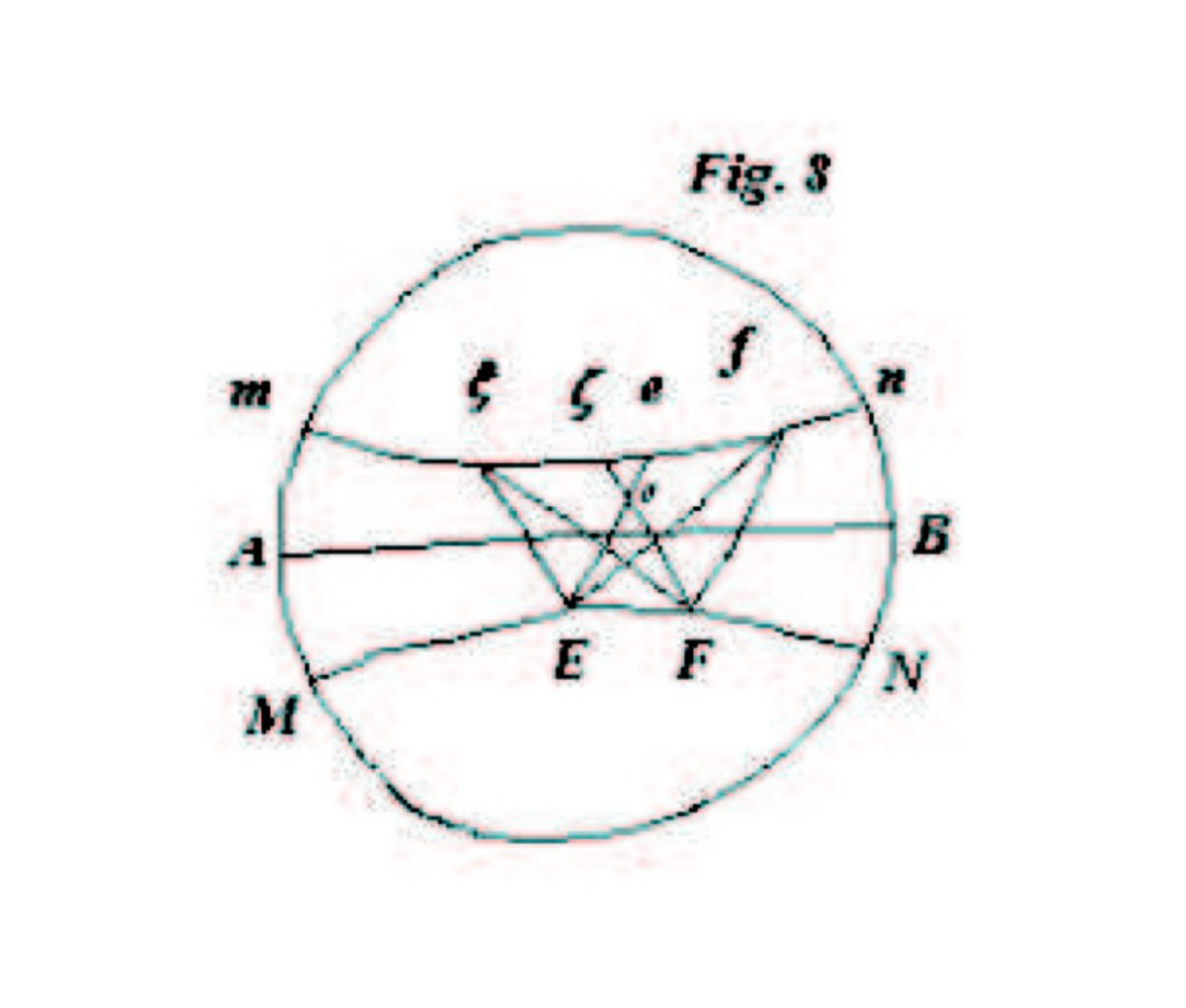}
\caption{\small {Extracted from Euler's paper \cite{Euler-Variae-T}.}}
\label{para1}
\end{figure}

 Euler's starts with a lemma. We use Euler's notation in Figure \ref{para1}, where $AB$ is an arc of great circle and $mn$ and $MN$ are two small circles which are equidistant of the great circle. The lemma says the following:
 \begin{quote}\small 
 If $eE$ is any arc joining the two small circles as in the figure  and if $O$ is its intersection with $AB$, then $eO=EO$ and $\widehat{meO}=\widehat{NEO}$.
 \end{quote}

We use again the notation of Figure \ref{para1}. Consider on one of our two small circles a  segment $EF$, and draw two segments of great circles $Ee$ and $Ff$, where $e$ and $f$ are on the other small circle, and making equal angles at $E$ and $F$.\footnote{Note that Euler uses here a notion of angle between arcs which are not necessarily arcs of great circles.} The Figure $EFfe$ is called a \emph{parallelogram} because any two opposite sides (even though they are not line segments of the geometry of the sphere) are equal, and any two opposite angles are also equal (Euler proves this). Taking now any other parallelogram $EF\zeta\xi$ having the same base $EF$, Euler shows that the two parallelograms $EFfe$ and $EF\zeta\xi$ have the same area.

 Each of the parallelograms $EFfe$ and $EF\zeta\xi$ breaks up into two triangles, by drawing the diagonals $Ef$ and $F\xi$. It follows from the preceding result that the four triangles $EfF$, $EeF$, $E\zeta F$ and $E\xi F$ have the area. Note that these may not be spherical triangles in the usual sense because some of their sides may not be arcs of great circles. However, we can enhance these triangles into genuine spherical triangles by replacing the side $EF$ by the segment of great circles that joins its extremities. Since each of these triangles is augmented by the same area, the four new genuine spherical triangles have the same area.
We deduce in particular the following:
\begin{quote}
\emph{The locus of the vertices of the spherical triangles whose common base is the arc of great circle $EF$ and which have the same area is the small circle $fe\zeta\xi$.}
\end{quote}

Having established this, Euler writes: ``By observing this, the problem due to the famous Professor Lexell may easily be resolved as follows" (\cite{Euler-Variae-T} \S 20), and he describes the following construction (We use Euler's notation in Figure \ref{Lexell2}):
 \begin{figure}[ht!]
\centering
 \psfrag{A}{\small $A$}
 \psfrag{B}{\small $B$}
 \psfrag{C}{\small $C$}
 \psfrag{p}{\small $p$}
 \psfrag{E}{\small $E$}
 \psfrag{P}{\small $P$}
 \psfrag{E}{\small $F$}
 \psfrag{E}{\small $e$}
 \psfrag{f}{\small $f$}
\includegraphics[width=.45\linewidth]{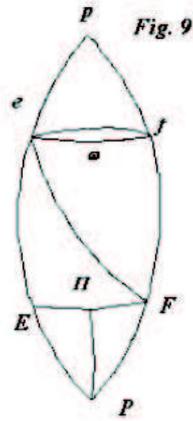}
\caption{\small {Extracted from Euler's paper \cite{Euler-Variae-T}}.}
\label{Lexell2}
\end{figure}
Let $EF$ be the base of the triangle that we seek for and let $\Delta$ be its area.
The vertex of the triangle varies on a small circle $ef$ whose center is a point $p$ and whose radius is $x$. Let $P$ be the point which is antipodal to $P$. We have $PE=PF=pe=pf$. Let $\phi=\widehat{PEF}$. Then we also have $\phi=\widehat{PFE}$, and Euler finds, for the area $\Delta$, the relation the relation $\Delta = \pi-2\phi$.  

Since $\Delta$ is given, we can find the point $P$ by erecting upon the base $EF$ an isosceles triangle having angles $\phi$ at the base. Therefore the point $P$ is known, and its antipode $p$ is also known.  The small circle $ef$ can then be drawn, since we know its center and its radius. 

Euler also notes that if we set $a=EF$, then the value $x$ of this radius is given by 
\[\tan x= \frac{\tan a/2}{\tan \Delta/2}.\]

 Let us now mention other developments of the same problem.
 
  In his paper \cite{Steiner}, published several decades after Lexell's paper, Jakob Steiner proves the following proposition:
\begin{quote}\small
The locus of the vertices of triangles of a given base and fixed area is a small circle that passes through the antipodes of the extremities of the base.
\end{quote}
The steps of his proof are as follows:

\begin{figure}[ht!]
\centering
 \psfrag{A}{\small $A$}
 \psfrag{B}{\small $B$}
 \psfrag{C}{\small $C$}
 \psfrag{D}{\small $D$}
 \psfrag{P}{\small $P$}
\includegraphics[width=.45\linewidth]{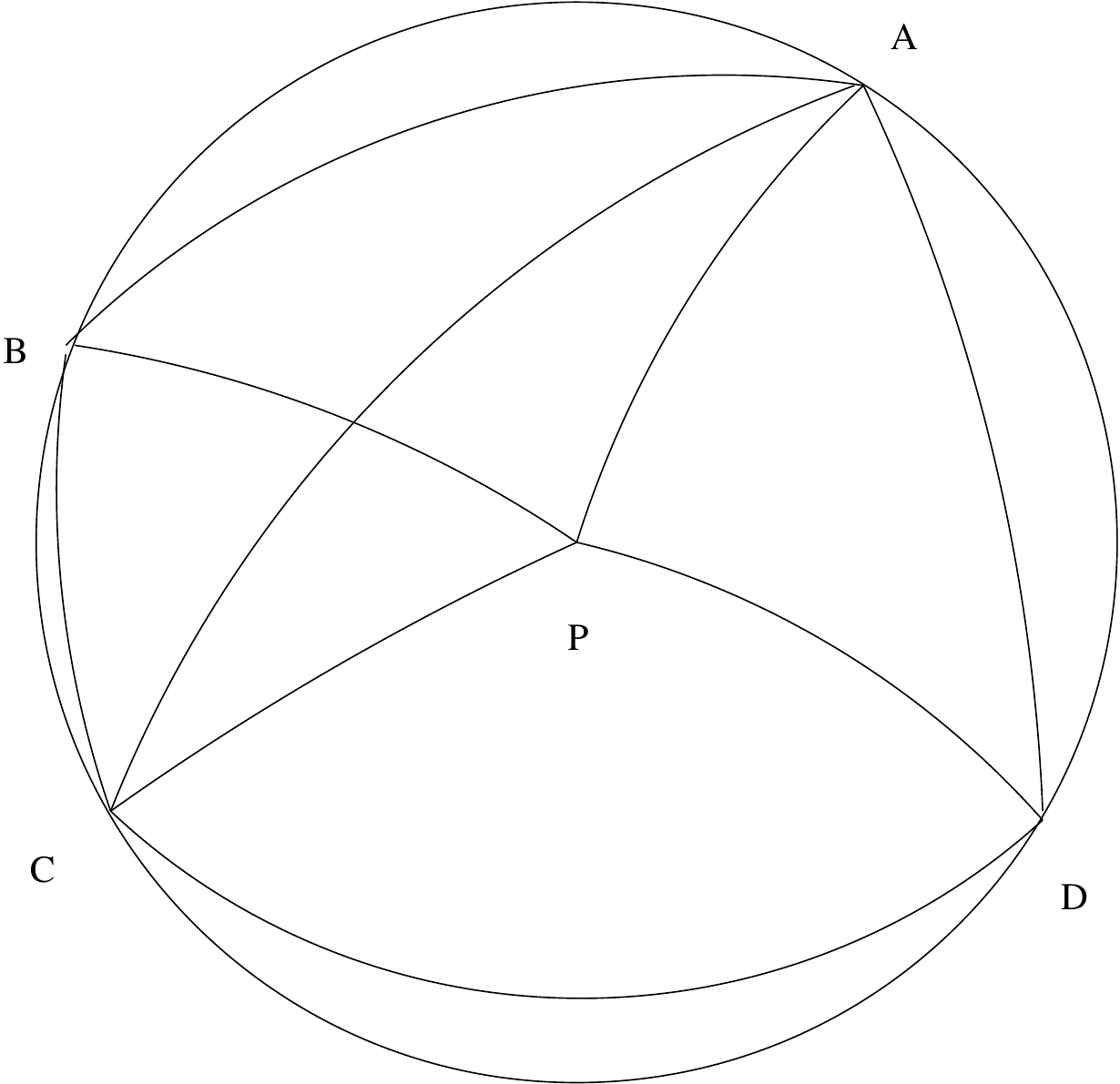}
\caption{\small {}}
\label{Steiner1}
\end{figure}
Consider an arbitrary quadrilateral $ABCD$ inscribed in a circle (Figure \ref{Steiner1}). Then, the two sums of opposite angles are equal, that is, $A+C=B+D$. This is seen by joining the pole $P$ of the circle $ABCD$ to each of the vertices $A,B,C,D$ and noting that the four triangles so obtained are isosceles and therefore their angles at the bases are equal.

Draw the diagonal $AC$ of this quadrilateral and let it cut the two angles at $A$ and $C$ into angles $\alpha,\alpha'$ and $\gamma,\gamma'$ respectively, such that the angles $\alpha$ and $\gamma$ belong to the triangle $BAC$ and $\alpha'$ and $\beta'$ belong to the triangle $DAC$. 
Then, we have 
\[\alpha+\alpha'+\gamma+\gamma'= B+D\]
and
\[(\alpha+\gamma)-B=D-(\alpha'+\gamma').\]

Keeping the three vertices $A,D,C$ fixed and varying $B$ on the arc $ABC$, the difference $\alpha+\gamma-B$, which is equal to $D-(\alpha'+\gamma')$ remains constant. Thus we have the following:
\begin{quote}
Given a spherical triangle $ABC$ with a fixed base $AC$ and with the difference between the angles at the base and at the vertex, $\alpha+\gamma-B$, fixed, the locus of the vertex is a circle $P$ passing through the extremities $A$ and $B$ of the base. 
\end{quote}
\begin{figure}[ht!]
\centering
 \psfrag{A}{\small $A$}
 \psfrag{B}{\small $B$}
 \psfrag{C}{\small $C$}
 \psfrag{A'}{\small $A'$}
 \psfrag{C'}{\small $C'$}
\includegraphics[width=.65\linewidth]{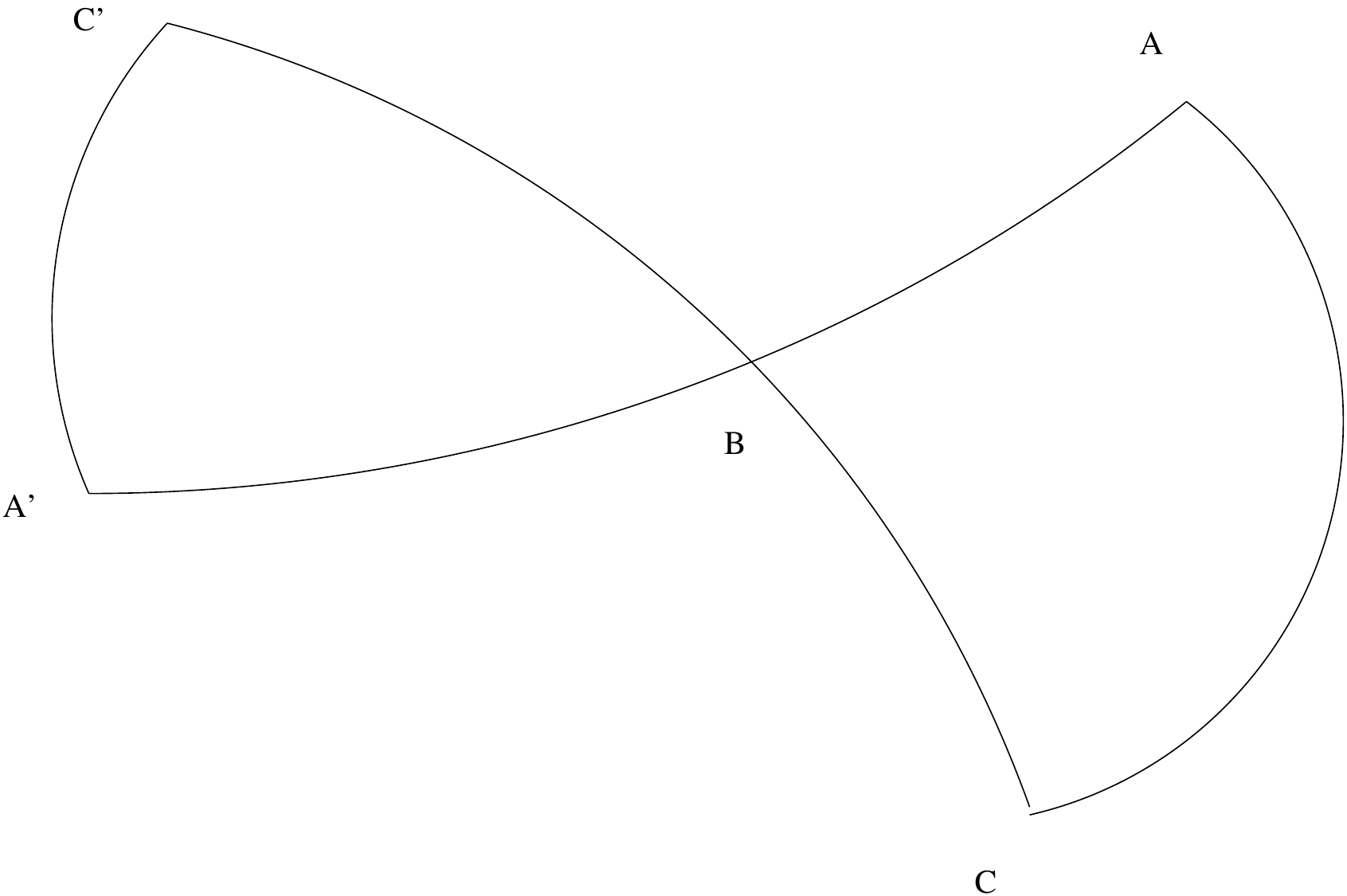}
\caption{\small {}}
\label{Steiner22}
\end{figure}
Extend the edges $AB$ and $CB$ of the triangle $ABC$ until the two points $A'$ and $B'$ which are the antipodes of $A$ and $B$ respectively, and consider the triangle $BA'C'$ (Figure \ref{Steiner22}). Then, we have
\[\alpha+A'=\pi=\gamma+C'.\]
Denoting by $K$ the constant $D-(\alpha'+\gamma')$, we get:
\[A'+B+C'=2\pi -K.\]
Thus, if the difference $\alpha+\gamma-B$ is constant in the triangle $ABC$, the angle sum, and therefore the area, of the triangle $A'BC'$ is also constant, and vice versa.
Thus, the locus of the point $B$ of the triangle $A'BC'$ is a circle which passes by the points $A$ and $C$.

Steiner writes in \cite{Steiner3} p. 101:
\begin{quote}\small
The history of this theorem presents a remarkable singularity. Due to Lexell, this theorem was generally known only through the \emph{\'El\'ements de G\'eom\'etrie} of Legendre who, although he attributed it to Lexell, gives it in an incomplete manner, and he seems to have been followed by all the authors who talked about it after him. Since I have been led in the memoir which I cited to recognize that ``the small circle which is the locus of the vertices of the triangles of equal area constructed over the same base passes through the two antipodes of the extremities of the base", I was led to think that this complement, which was essential for the applications which I had in mind, was not known, and I was comforted in this error by all those who later on dealt with the same subject. It is only recently that Mr. Liouville, who made a review of the present Memoir to the Academy of Sciences of Paris, and having had the idea of resorting to the original memoir of Lexell (Nova Acta Petrop. Tom. V Pars Prima), recognized that the proposition in question is stated there in a complete manner and is proved in two different ways. It is not possible to guess why Legendre mutilated the theorem given by Lexell, and one is the more surprised that another circumstance remained unseen for a long time, namely, that the same proposition was the subject of a memoir of Euler (Nova Acta Tom. X.) where it is proved in a very elegant and purely geometric manner. I would add that the proof given by this famous geometer presents several analogies with the one I indicated in the first publication of the present memoir, in the journal founded by Mr. Liouville and which is based on considerations pertaining to three-dimensional geometry.
\end{quote}

In his paper \cite{Lebesgue1855}, Lebesgue gives a proof of Lexell's theorem which in fact is Euler's proof. At the end of this paper, the editor of the journal adds a comment, saying that one can find a proof of this theorem in the \emph{\'El\'ements de G\'eom\'etrie} of Catalan (Book VII, Problem VII), but, again,  no reference is given to Euler. 

The only original contribution of Lebesgue in this paper consists in two remarks that he makes after the proof. The first remark says that one has an analogous result for surfaces of revolution which admit an equator, taking two curves that are equidistant  to that equator. The second remark says that by taking the stereographic projection of the figure on the sphere, we get in the plane a family of curved triangles on a given base, whose angle sum is constant, and whose vertices are on the same circle.

A solution of the question underlying Lexell's theorem can also be given by exhibiting an explicit formula for the area of a triangle of a fixed base and a variable third vertex. Such a formula leads to an analytic equation for the third vertex and one can eventually show that this is the equation of a circle on the sphere. Euler and Lexell also followed this path in their papers  \cite{Euler-Variae-T} and  \cite{Lexell-Solutio}. Another such approach  was conducted by Puissant in his \emph{Trait\'e de  g\'eod\'esie, ou exposition des m\'ethodes trigonom\'etriques et astronomiques} \cite{Puissant} (p. 115).

\section{Area}\label{s:area}

In this section, we present some important formulae by Euler and his followers for the areas of spherical triangles. 

The formula for areas of triangles is attributed to Albert Girard (1595-1632), who stated it in his \emph{Invention nouvelle en alg\`ebre}  (1629). Lagrange says in his paper \cite{Lagrange-Solution-T} that the solution of Girard is not satisfying, and that the theorem should rather be attributed to Cavalieri, who gave a compete proof in his \emph{Directorium generale uranometricum} (A general guide to celestial measures, Bologna, 1632). Lagrange, in the same paper, also refers to the ``beautiful proof" of Wallis (1616-1703). 

The formula for the area of a triangle in terms of angular excess was also found by Thomas Harriot (1560-1621) who did not publish it.

\begin{figure}[ht!]
\centering
 \includegraphics[width=.90\linewidth]{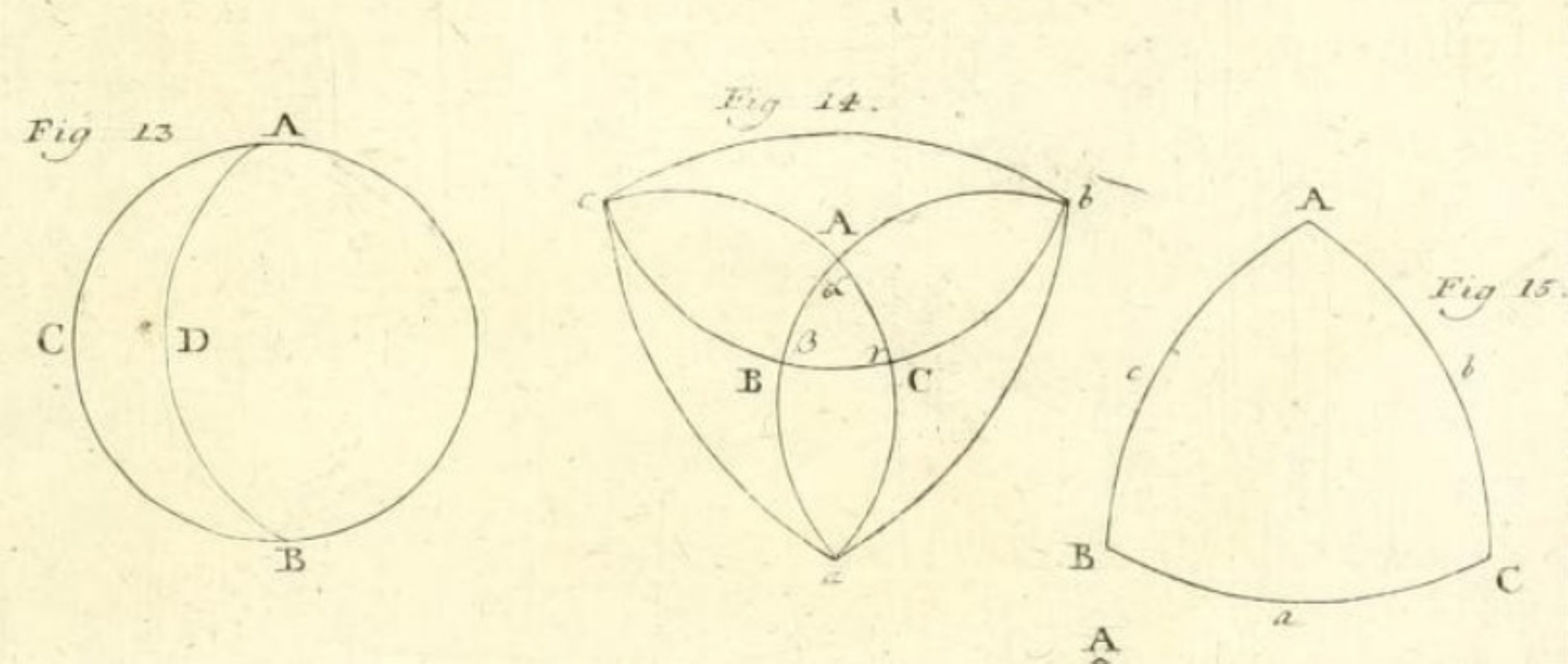}
\caption{\small {Extracted from Euler's original memoir \cite{Euler-Mensura-T} in which Euler proves Girard's theorem using the notion of polar triangle, drawn in the middle of this plate.}}
\label{polar1}
\end{figure}

In his memoir
\emph{Principes de la trigonom\'etrie sph\'erique tir\'es de la m\'ethode des plus grands et des plus petits} \cite{Euler-Principes-T} which we already mentioned, Euler gives a proof of this theorem using the calculus of variations. In his memoir \emph{De mensura angulorum solidorum} (On the measure of solid angles) \cite{Euler-Mensura-T}, written in 1775, he gives another proof of this theorem based on the notion of polar triangle\footnote{We recal that the polar triangle of a triangle $ABC$ is a spherical triangle $A'B'C'$ whose sides $B'C', C'A', A'B'$ are on arcs of great circles with poles $A,B,C$. (Thus, $A$ is a pole of the equator $B'C'$, etc.)}  (see Figure \ref{polar1}).  He starts with the fact that the area of a lunar part\footnote{A \emph{lune} is a region such as $ACBDA$ in Figure \ref{polar1}.} of the sphere (a region comprised between two arcs of great circle joining the same
 points) of angle $A$ is equal to $2A$. The formula for the area of a triangle follows then from drawing the polar triangle, adding the areas of the three lunes obtained (whose angles are those of the given triangle) and subtracting the value of four areas of spherical triangles which can be easily compared to the area of the initial triangle. Lagrange, in his 1800 memoir, does not mention the proof of Euler.
We also mention that the existence of the polar triangle is at the basis of the symmetry in the trigonometric formulae which is manifested by the fact that some of them can be obtained from  others by exchanging the roles of angles and sides, as we already mentioned in \S \ref{sph}.

Let us now elaborate on the formulae for area in terms of side lengths.

In Euclidean geometry, the formula which gives the area of a triangle in term of the length of its sides is attributed to Heron (c. 10--70 A.D.), the Greek mathematician who worked, like Euclid, in Alexandria.

 \begin{figure}[ht!]
\centering
 \psfrag{A}{\small $A$}
 \psfrag{B}{\small $B$}
 \psfrag{C}{\small $C$}
 \psfrag{R}{\small $R$}
 \psfrag{V}{\small $V$}
\includegraphics[width=.40\linewidth]{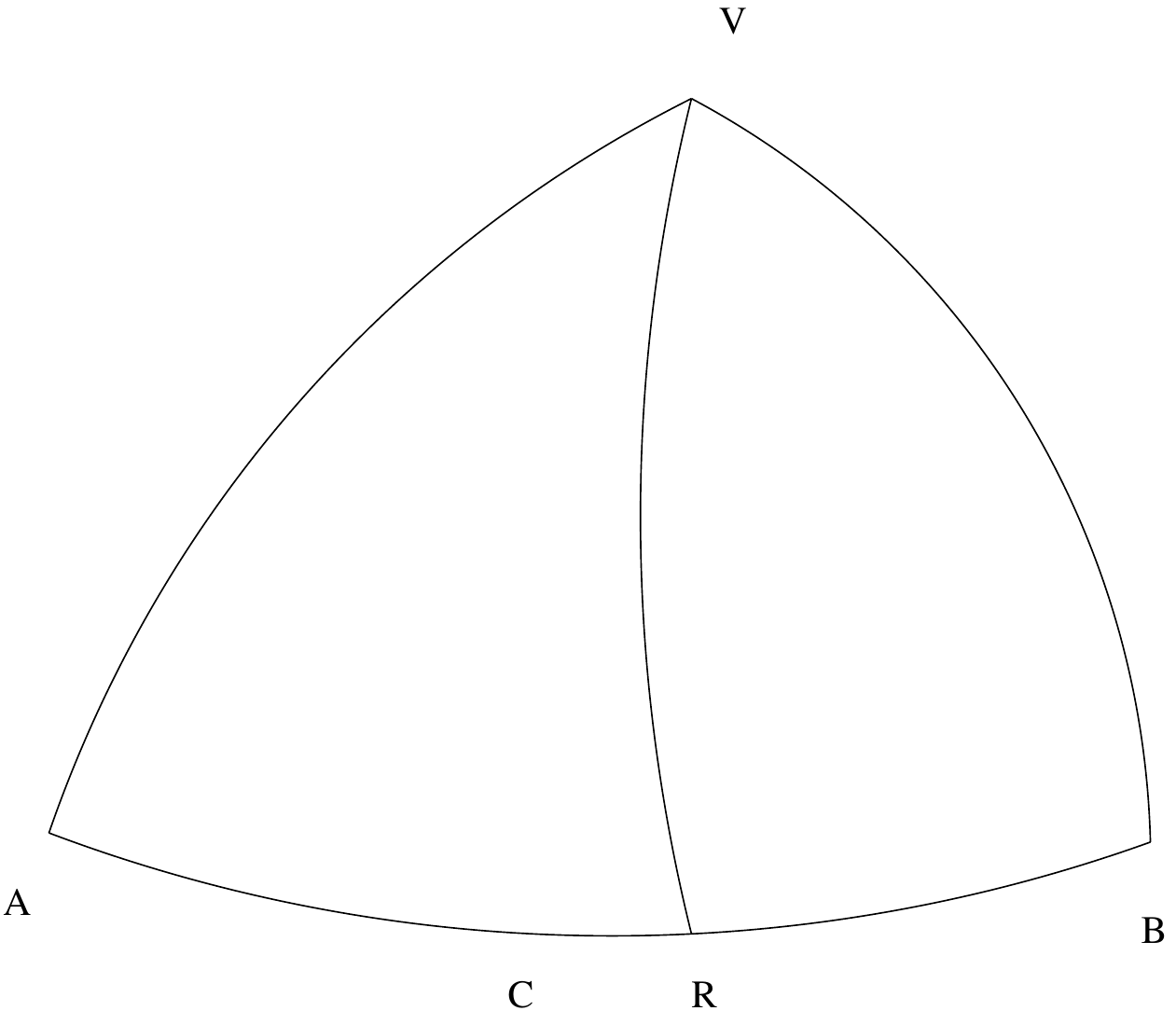}
\caption{\small {}}
\label{Area-Lexell}
\end{figure}

Using Lexell's notation  \cite{Lexell-Solutio} (Figure \ref{Area-Lexell} here), we consider a triangle $ABV$ with sides $a,b,v$ opposite to the vertices $A,B,V$ respectively.
We let $C$ be the midpoint of $AB$ and we set $VR=y$ and $CR=x$.

Lexell gives the following formula for the area $\delta$ of $ABV$:  
\[\cot \delta = \frac{\cos y \cos x + \cos a}{\sin a \sin y}.\]

In his memoir \cite{Euler-Mensura-T}, Euler gives several formulae for the area in terms of side lengths. Denoting the sides of the triangle by $a,b,c$, he gives the following formula: 
\[
\tan \frac{1}{2} \delta=\frac{\sqrt{1-\cos^2a-\cos^2b-\cos^2c+2\cos a\cos b\cos c}}{1+\cos a+\cos b+\cos c}.
\]
In \cite{Euler-Variae-T}, he gives the following formula:
\begin{equation} \label{eq:area:Euler}
\cos \frac{1}{2}\delta = \frac{1+\cos a +\cos b +\cos c}{4 \cos\frac{1}{2} a \cos\frac{1}{2} b\cos\frac{1}{2} c}.
\end{equation}
The paper is written in 1778 and Euler says that he found the formula a long time before. He considers this as an outstanding theorem (egregium theorema).

\begin{figure}[ht!]
\centering
 \psfrag{A}{\small $A$}
 \psfrag{B}{\small $B$}
 \psfrag{Z}{\small $Z$}
 \psfrag{z}{\small $z$}
 \psfrag{V}{\small $V$}
\includegraphics[width=.55\linewidth]{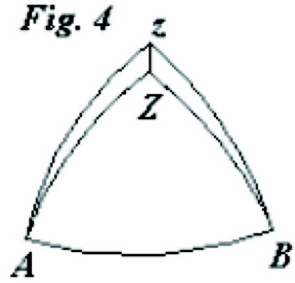}
\caption{\small {Extracted from Euler's \emph{Variae speculationes super area triangulorum sphaericorum}. The vertex $Z$ undergoes an infinitesimal variation.}}
\label{Euler-Variae1}
\end{figure}

 Euler gives two proofs of Formula (\ref{eq:area:Euler}). The first one uses analysis (differential calculus) and it is in the spirit of the methods of the calculus of variations. Euler starts by taking a spherical triangle $AZB$ with base $AB$ and varying the lengths of the sides $AZ=x$ and $BZ=y$ by an infinitesimal amounts, thus getting a new triangle $AzB$ (see Figure \ref{Euler-Variae1})
 whose sides are $Az=x+\partial x$ and $Bz=y+\partial y$. He then investigates the area of this new triangle. Denoting by $\Delta$ the area of this triangle and $\phi$ and $\psi$ the angles $\widehat{BAZ}$ and $\widehat{ABZ}$ respectively, he obtains the equation
\[\partial \Delta = \partial \phi (1-\cos x) +\partial \psi (1-\cos y).\]
He then eliminates the angles $\phi$ and $\psi$ by using the trigonometric formulae, he gets a formula for $\partial \Delta$ which only involves side lengths. After several integrations, he gets Formula (\ref{eq:area:Euler}).

 The second proof by Euler of the area formula (\ref{eq:area:Euler}) uses  Girard's Theorem which we mentioned, that is, the formula which gives the area as the angular excess. In particular, combining the two proofs, we get a new proof of the theorem of Girard. Several other formulae for the area of a triangle in terms of the side lengths were obtained by Euler.

 Denoting, as before, the area by $\delta$ of a spherical triangle of sides $a,b,c$, contained in a paper by Lagrange which was published in 1800 (\cite{Lagrange-Solution-T}, p. 340):

\[\tan \frac{\delta}{2}=\frac{{2\sqrt{\sin(\frac{a+b+c}{2})\sin(\frac{a-b+c}{2})\sin(\frac{a+b-c}{2})\sin(\frac{-a+b+c}{2})}}}{1+\cos a+\cos b+\cos c}.\]

 In his \emph{Trait\'e de  g\'eod\'esie} \cite{Puissant}, Puissant gives several expressions of the area of a spherical triangle, for example the following, which gives the value of the area  in terms of two sides and the angle they contain (p. 109):
\[\tan\frac{1}{2}\delta = \frac{\tan \frac{1}{2}a\tan \frac{1}{2}b\sin C}{1+\tan \frac{1}{2}a\tan \frac{1}{2}b \cos C}.
\]
and then (p. 110)
\[\cot \frac{1}{2}\delta = \frac{1+\cos a +\cos b+\cos c}{\sin a\sin b\sin C}.\]
These formulae have interesting applications.

\section{On a problem of Pappus}\label{s:pappus}

We already stated this problem in \S \ref{s:quick} . It asks for the construction of a triangle inscribed in a circle whose sides are contained in three lines that pass through three given points.
\begin{figure}[ht!]
\centering
 \includegraphics[width=.90\linewidth]{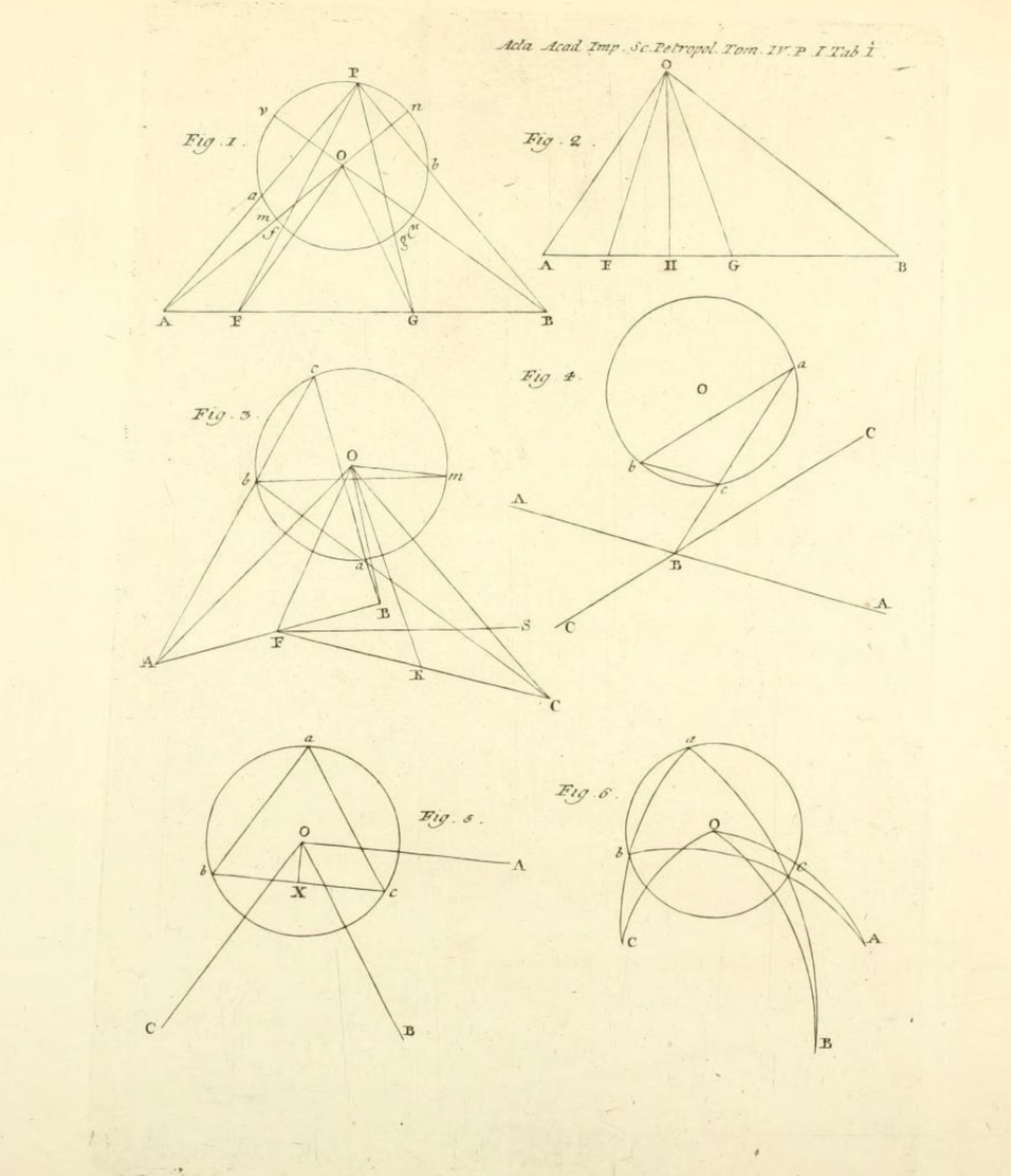}
\caption{\small {Plate extracted from Euler's original memoir \cite{Euler-Pappi-T} with the figures corresponding to the various versions of Pappus' construction and its generalizations. The last one is the spherical case.}}
\label{Pappus}
\end{figure}
     
  The problem has several ramifications. In the same volume of the  \emph{Acta academiae scientarum imperialis Petropolitinae} in which Euler's soultion appeared, another solution of the problem of Pappus was given by his student Fuss \cite{Fuss-Solutio}.
  
  Several other major mathematicians worked on the same problem. Chasles, in his \emph{Aper\c cu historique} \cite{Chasles-Apercu} p. 328 has a note on the history of this problem (in the Euclidean setting). We learn from him that in 1742, Cramer proposed the general problem (in the Euclidean case) to Castillon and the latter published a solution, in 1776, in the Berlin memoirs. The same volume contains a purely analytic solution by Lagrange (using trigonometric formulae).  Giordano di Oltaiano and Francisco Malfatti solved the more general question for a polygon with an arbitrary number of sides that pass by an equal number of points. Their works are published in Volume IV of the \emph{Memorie della societa italiana}.  In the Berlin memoirs of 1796, l'Huilier  published a modified algebraic solution of the same general problem. Carnot (p. 383 of his \emph{G\'eom\'etrie de position}) gave a modified form of Lagrange's solution. He applied Lagrange's solution to the case of an arbitrary polygon.  Brianchon, in the Journal de l'\'Ecole Polytechnique, 10e cahier, published a solution in the case where a conic section replaces the circle, and where  the three points are on a line. Gergonne worked out a solution in the case where the circle is replaced by a general conic and where the three points are not necessarily aligned, and he allows only the straightedge (and no compass) for the solution (\emph{Annales de Math\'ematiques}, tom. 1, p. 341, 1810-1811). In fact, Gergonne first considered the problem of circumscribing in a conic a triangle whose vertices are on three given lines, and he later on observed that using polarity theory for conics, one can transform the problem into the one of circumscribing in a conic a triangle whose vertices pass by three given points. Gergonne considers this problem, together with other classical problems (from Pappus and others)
  in the setting of analytic geometry. An analysis of his work in this matter is contained in \cite{Notice} p. 168ff. Finally, Poncelet, in his \emph{Trait\'e des propri\'et\'es projectives} (p. 352) considers the extension to the case of a polygon with an arbitrary number of sides, and where the circle is replaced by an arbitrary conic.

\section{Sphere projection and cartography}\label{s:maps}

Spherical geometry is a mathematical theory with applications, and Euler was interested in these applications. Chasles writes in his  \emph{Aper\c cu historique} \cite{Chasles-Apercu} (p. 235) that since the works of the ancients (he mentions Theodosius, Menelaus and Ptolemy), ``if the theory [of spherical geometry] , was extended and attained, in the hands of our most celebrated geometers, a high degree of perfection, it was always done by almost preserving the same framework, because the goal was always the same: the computation of triangles in order to fit the service of the astronomer and the navigator, and for the great geodesic operations which showed the true form of the Earth."

In the year 1777, Euler wrote the three memoirs on maps from the sphere to the Euclidean plane \cite{Euler-rep-1777},  \cite{Euler-pro-1777} and   \cite{Euler-pro-Desli-1777}.\footnote{The title of the last memoir,  \emph{De proiectione geographica Deslisliana in mappa generali imperii russici usitata},  refers to Joseph-Nicolas Delisle (1688-1768), a leading French astronomer and geographer who worked at the Saint Petersburg Academy of Sciences. Having been called by the emperor Peter the Great  at the founding of the Academy, he arrived there in 1726 --  one year before Euler. Delisle founded and run the observatory of Saint Petersburg, on the Vasilyevsky Island, which was one of the finest in Europe, and he was also in charge of composing maps of Russia. Between 1726 and 1747, Delisle was a member of the Saint Petersburg Academy of Sciences, and at the beginning of his stay there, Euler assisted him in recording astronomical observations which were used to make meridian tables. Delisle returned to Paris in 1747, where he founded the famous observatory at the h\^otel de Cluny.}
. The three papers are motivated by the practical question of drawing geographical maps. We may recall at this point that Euler, at the Academy of Sciences of Saint Petersburg, besides being a mathematician, also held the official charge of cartographer. We also recall that knowing the exact form of the Earth was one of the major concerns of the eighteenth-century scientists, and there were long debates and controversies on this matter opposing the French and the English scientists. This also led to famous expeditions, to Peru, Lapland, Kamchatka and elsewhere. It is also important to recall that scientists in Europe  were given by the great monarchs (Frederick II, Louis XIV and the Russian tsars)  the task of finding methods in order to measure exactly the size of their kingdoms.  Several methods were developed for that purpose, and they used triangulations.

The paper \cite{Euler-rep-1777} is the longest among the three. Euler studies there systematically several maps from the sphere to the Euclidean plane, having different desired properties. 
The method is through differential calculus, that is, he finds differential equations which the required maps satisfy.

Euler starts by showing a fact which he considers well-known, which is that there is no ``perfect" or ``exact" mapping from the sphere onto the plane. In modern words, Euler shows that the sphere is not developable (even locally) onto the plane.  Having proved this fact, he writes:\footnote{For this paper, we use W. Heine's translation.} ``Since therefore a perfectly exact representation is excluded, we are obliged to admit representations which are not similar, so that the spherical figure differs in some manner from its image in the plane."

Euler studies the  three sorts of maps:

\begin{enumerate}
\item \label{map1} Maps where images of all meridians are normal to a given axis (the ``horizontal" axis in the plane), while all parallels are sent parallel to it. 

\item \label{map3} The maps which ``preserves the properties of small figures", that is, that areconformal.

\item \label{map2} Maps where surface area is represented at its true size.

\end{enumerate}

Examples of maps satisfying Property (\ref{map1}) are the famous Mercator maritime charts, which Euler mentions, and which are obtained by projecting the sphere onto a cylinder which is tangent to the equator. Distances can be highly deformed by such maps (at the poles, they become infinite), but there are certain practical advantages, and these maps are still used in sailing, because straight lines drawn on the map represent lines of constant slope on the globe.\footnote{That is, the angle these lines makes with the meridians is constant.} Sailors are used to follow these paths, even if they are not geodesics.\footnote{The lines are practical to follow on a Mercator maritime map, because they are straight lines. At the poles, these lines are spirals and have infinite length. (This is the meaning of the word ``loxodromic" which Euler uses in the next sentence.) The use of the Mercator maps is limited to latitudes which are below $70^{\mathrm{o}}$.} Euler  writes that ``the greatest advantage which this map gives to travelers at sea is that the loxodromic curves, which on the sphere cut each meridian at the same angle, are here represented by a straight line. Such a straight line cuts all the meridians, which are parallel to each other, at the same angle [...] But while loxodromoic curves are represented on the plane simply as straight lines, in contrast, great circles on the sphere are represented by transcendental curves of a very high level."

 Euler concludes his paper with the following: ``In these three hypotheses\footnote{This refers to Properties (1) to (3) that we mantioned above.} is contained everything ordinarily desired from geographic as well as hydrographic maps. The second hypothesis treated above even covers all possible representations. But on account of the great generality of the resulting formulae, it is not easy to elicit from them any methods of practical use. Nor, indeed, was the intention of the present work to go into practical uses, especially since, with the usual projections, these matters have been explained in detail by others."

In his memoirs \cite{S0}, \cite{S01} and \cite{S00}, Schubert, another student of Euler, examines the errors which one encounters when considering the terrestrial globe as a sphere. He also describes the corrections that have to be made if one considers the globe as a rotationally symmetric spheroid.\footnote{This is a surface obtained by the rotation of an ellipse around one of its axes. Not all spheriods have such a symmetry.} The projections of the meridians become ellipses instead of circles (as they are in the spherical case). Lambert gives precise formulae for the axes of these spheroids and for their intersections with the meridian. He also studies the question of mapping conformally the spheroid, that is, in an angle-preserving way, onto a plane. In the  memoir \cite{S00}, he examines geographical maps obtained from maps between a sphere and the surface of a cone.

Mathematical questions related to the problem of drawing geographical maps were thoroughly  studied by Lambert,  Lagrange, Gauss, Beltrami, Darboux, Liouville, Bonnet and others.  Lagrange wrote two memoirs on this subject \cite{Lagrange1} extending works on Lambert and Euler. He considers (cf. the Introduction to his First Memoir) that Lambert was the first to consider the theory of geographic maps in full generality. Lambert\footnote{The Alsatian mathematician Johann Heinrich Lambert  (1728-1777) is sometimes considered as the founder of modern cartography.  His \emph{Anmerkungen und Zus\"atze zur Entwerfung der Land- und Himmelscharten} (Remarks and complements for the design of terrestrial and celestial maps, 1772) \cite{Lamb-Anmer} contains seven new map projections, including the Lambert conformal conic projection, the transverse Mercator, the Lambert azimuthal equal area projection, and the Lambert cylindrical equal-area projection. Several among these projections are still in use today, for various purposes. Lambert was one of the most brilliant precursors of hyperbolic geometry. In his {\it Theorie der Parallellinien}, written in 1766, he developed the bases of a geometry in which all the Euclidean postulates hold except the parallel postulate which is replaced by its negation. Lambert's work is edited in the collection  \cite{Lambert-Para}. We refer the reader to \cite{Lambert-Blanchard} for the first translation of this work, together with a mathematical commentary. Euler had a great respect for Lambert, and he helped him joining the Academy of Sciences of Berlin.}, in his famous paper \emph{Beitr\"age zum Gebrauche der Mathematik und deren Anwendung} (Contributions to the use of mathematics and its applications) \cite{Lambert-Bey} studied the question of characterizing the projections of the sphere onto the plane such that the angles between the images of the meridians and the parallels are preserved. The solution given  by Euler to that problem followes that of Lambert.
In 1816, H. C. Schumacher, a famous German-Danish astronomer at Copenhagen, who was a friend of Gauss (and who had been his student in G\"ottingen), announced a prize, to be awarded by the Royal Society of Copenhagen, for the question of obtaining a general conformal map between two surfaces. Gauss won the prize, for a solution he submitted in 1822, which was publised in 1923 in the \emph{Astronomische Nachrichten}. There exists an easily available French translation of this memoir of Gauss \cite{Gauss-s}. 
Beltrami worked on such questions  in the setting of the differential geometry of surfaces. In the introduction to his paper \cite{Beltrami1865}, he writes that a large part of the research done before him on such questions was concerned with conservation either of angles or of area. He adds that even though these two properties are regarded as the simplest and most important properties for geographic maps, there are other properties that one might want to preserve. He declares that since the projection maps that are used in this science are  mainly concerned with the measurement of distances, one would like to exclude projection maps where the images of distance-minimizing curves are too remote from straight lines on the plane. He mentions in passing that the central projection of the sphere is the only map that transforms the sphere geodesics into Euclidean straight lines.  Beltrami then says that beyond its applications to geographic map drawing, the solution of the problem may lead to ``a new method of geodesic calculus, in which the questions concerning geodesic triangles on surfaces can all be reduced to simple questions of plane trigonometry." 

 \section{On curves on the surface of the sphere}\label{s:alg}
 
We mention briefly the memoir \cite{Euler-Curva1771} in which Euler studies the existence of curves on the sphere that can be expressible by algebraic equations.\footnote{This is the meaning of the word ``rectifiable" in the title of the paper. Note that this is different from the sense that this word has today in mathematics.} There are two other papers on the same subject, \cite{Euler-coni} and \cite{Euler-Lineis} where the same kind of questions are studied on other surfaces, namely, the cone and the spheroid. In the paper  \cite{Euler-coni}, Euler finds a family of rectifiable curves on a right cone whose ratio of side-length to the diameter of the base is rational. In the paper \cite{Euler-Curva1771}, he finds such a curve on the sphere, namely, the ``spherical epicycloid", obtained by the motion of a great circle on a small circle whose diameter makes a rational ratio with the diameter of the sphere. The paper \cite{Euler-Lineis}, on the spheroid, gives a simpler proof of the result of \cite{Euler-Curva1771} and it generalizes it to the case of curves on an arbitrary spheroid.

\section{Notes on Euler's students}\label{s:students}

This last section contains some information on Euler's students and young collaborators, especially those who continued his work on spherical geometry.

At the two Academies to which he belonged, Euler never had an official teaching duty, and he only occasionally gave lessons, notably to princes and to princesses, in particular at the beginning of his stay at Berlin  (cf.  
 \cite{Fuss-Eloge} note (f) p.~48). 
Euler  nevertheless had evident pedagogical abilities, as one can see by skimming in the books and treatises he wrote. He even wrote school textbooks. We refer the reader to Chapter 1 of \cite{HP}.

We learn from Fuss' \emph{\'Eloge fun\`ebre}  \cite[p. 72]{Fuss-Eloge},  that at the time of Euler's death, there were, at the Academy of Saint Petersburg,  eight mathematicians who had benefited from his teaching, namely, his oldest son Johann Albrecht  Euler\footnote{\label{f:JA}Leonhard Euler had thirteen children, but only five of them attained adult age.  Johann Albrecht (1734-1800) was his oldest son.  He was an excellent scientist, although his career evolved in the shade of his father. He was elected, in 1754, at the age of 20, at the Academy of Sciences of Berlin, of which his father was member. In 1765, when the Euler family returned to Russia, Johann Albrecht became member of the Academy of Sciences of Saint Petersburg, at the chair of physics. Together with Lexell and Fuss, he helped his father during the long period where the latter was blind, assisting him in preparing his papers, writing them up and reading them at the Academy. In Saint Petersburg, the family of Johann Albrecht lived in the house of the father, occupying the base floor.}, Semion Kirillovich Kotelnikov, Stepan Rumovsky, Georg Wolfgang Krafft, Anders Johan Lexell, Petr  Inokhodtsov, Mikhail Evseyevich Golovine and Nicolaus Fuss, the author of the \emph{\'Eloge}. These eight scientists are sometimes considered to be the students of Euler. They learned from him, but they also assisted him in preparing his papers and translating his books into Russian. For instance, Inokhodtsov translated Euler's \emph{Elements of Algebra}, Rumovsky translated his \emph{Letters to a German Princess},  Golovine translated his \emph{Theorie complette de la construction et de la man\oe uvre des vaisseaux}. Krafft helped Euler writing his three-volume \emph{Dioptrica} and Fuss helped Euler with the preparation of over two hundred and fifty works. We also mention that Kotelnikov led a group that translated Christian Wolff's multi-volume \emph{Anfangsgr\"unde aller mathematischen Wissenschaften}.

Concerning the number of Euler's students, Condorcet is more generous. He declares, in his  \emph{\'Eloge} (\cite{Condorcet-Eloge},  p. 308): 
``All the famous mathematicians that exist today are his students: there is no one among them who has not been shaped by reading his books, or by receiving from him the formulae or the method he used, or who, in his discoveries, has not been guided by the spirit of Euler." From Condorcet, we also learn that ``among the sixteen professors that were attached at that time to the Academy of Sciences of Saint Petersburg, eight were shaped by Euler, and all of them, who were known from their works and who had received awards and academic titles, glorified in being Euler's disciples."

Some of Euler's students taught at various Russian schools and educational institutions, transmitting the teaching of their master to the future Russian generations.

At the end of his life, Euler left about 200 unfinished memoirs, and some of them were revised by his students.

Some details on the achievements of Euler's students are contained in Vucinich \cite{Vucinich-Mathematics-T} \cite{Vucinich-T} and in Fellmann \cite{Fellmann2}.   Five among these students -- Fuss,  Lexell, Schubert, Rumovsky and Golovine -- worked on spherical geometry, developing the ideas of their master. We already mentioned works of the first three. Rumovsky became an astronomer and a mathematician, and he headed the Saint Petersburg Geography department and astronomical observatory. Golovine published in 1789 a book on \emph{Plane and spherical geometry with algebraic proofs} which remained during several years the main reference on the subject. Rumovski became a renown astronomer and he was elected at the Swedish Academy of Sciences. He is the author of very precise geographical tables of the Russian Empire.  

We also learn from Fuss' Eulogy that Euler liked to work with his students, and he invited them regularly at his home, to discuss various subjects around a dinner table. During the years 1752-1756, while he was in Berlin, two Russian students followed him there, Kotelnikov and Rumovsky, and they were lodged at his house.

 Paul Heinrich Fuss\footnote{ Paul Heinrich Fuss (1797-1855) was a mathematician and physicist, the son of Nicolaus Fuss and of Euler's grand-daughter Albertine. He became adjunct at the Academy of Sciences of Saint-Petersburg in 1818, and full member in 1823. He became permanent secretary of the Academy in 1826. He published  the scientific part of the correspondence between  Euler, Johann Bernoulli, Nicolaus Bernoulli, Daniel Bernoulli and Christian Goldbach,  and  a short biography of Euler and and bibliography of his works, see \cite{Fuss-Corresp}.},  writes in his \emph{Correspondance} \cite{Fuss-Corresp} p. {\sc xliv}.\footnote{My translation.} ``It seems that before the arrival of my father, which took place in May 1773, Euler was helped some times by one, and some times by another one of his numerous students which were among his colleagues. In a big in-folio which I keep cherishingly and which contains the first sketches of Euler's memoirs which are prior to the mentioned epoch, I think I can recognize above all the hand of Krafft, as well as the ones of J.-A. Euler and of Lexell.  But I also notice that   they often contented themselves of simply sketching his memoirs, without going into a lot of trouble  to finish the writing.  [...] Every morning, a student showed up to read for him, either his vast correspondence (and he was totally in charge of it), or of political sheets, or of some work which was worthy of interest; they discussed several subjects of science, and the master, at this occasion, submitted with good grace to clarify things and to solve the difficulties of matters which the student had encountered during his studies. When the table was full of calculations -- which was often the case -- the master confided to his disciple his recent and fresh conceptions, and he showed him the course of his ideas and the general plan of the writing, leaving for him the development of the calculations; and usually, the student would bring to him  the following day the sketch of the memoir which was registered in the big book which we already mentioned (\emph{Adversaria mathematica}). Once this sketch was approved, the piece was written up neatly and was immediately presented to the Academy."

We finally mention some details on the lives and works of Fuss, Lexell and Schubert. The  three were mentioned several times in the previous sections.

Nikolaus Fuss (1755-1826) was Swiss, born in Basel where he received his education. Euler, who was nearly blind and who needed a secretary, asked his old friend Daniel Bernoulli to send him one from his hometown, and Bernoulli recommended Fuss, who arrived to Saint Petersburg in 1773. In 1776, Fuss was appointed junior scientist assistant at the Academy of Saint Petersburg. In 1783, at the death of Euler, he became an academician in the class of mathematics. 
From 1800 and until his death in 1826, he was the permanent secretary of the Academy.  During the seven years where he assisted Euler, Fuss helped him in the preparation of more than 250 memoirs. Fuss married in 1784 the grand-daughter of Euler (Albertine, the oldest daughter of his son Johann Albrecht). He is considered to have been the favorite student of Euler and a large part of his mathematical results are solution of problems raised by his master. His work concerns the theory of series, curves, differential equations, mechanics, astronomy and geodesy, and above all geometry.

Fuss also became honorary member of the academies of Berlin, Sweden and Danemark. He taught at the military academy, and (like Lexell), he taught at the naval cadet academy. He wrote several textbooks on elementary mathematics and he contributed to raising the Russian academic system. 

Certain works of Euler that were prepared by Fuss may be attributed to both of them.  The titles of the following papers by Fuss are significant in this respect: 

$\bullet$  \emph{ Instruction d\'etaill\'ee pour porter les lunettes de toutes les diff\'erentes esp\`eces au plus haut degr\'e de perfection dont elles sont susceptibles,  tir\'ee de la th\'eorie dioptrique de M. Euler le p\`ere et mise \`a la port\'ee de tous les ouvriers en ce genre par Nicolaus Fuss. Avec la description d'un microscope qui peut passer pour le plus parfait dans son esp\`ece et qui est propre \`a produire tous les grossissements qu'on voudra} (Detailed instructions in order to lead glasses of all different species to the highest possible degree of perfection, extracted from the theory of dioptrics of Mr. Euler the father, and  put within the reach of all the workers in this domain by Nicolaus Fuss. With the description of a microscope that can be considered as the most perfect of this kind and which is appropriate for producing all the desired magnifications).  Saint Petersburg, 1774 

$\bullet$ \emph{\'Eclaircissemens sur les \'etablissemens publics en faveur tant des veuves que des morts avec la description d'une nouvelle esp\`ece de tontine aussi favorable au public qu'utile \`a l'\'etat calcul\'es sous la direction de Monsieur Leonard Euler} (Enlightenment for the public institutions for the sake of supporting the widows as well as  the dead, with a description of a new species of tontine which is as much favorable to the public as useful to the state, calculated under the direction of Mr. Leonard Euler).  Saint Petersburg, 1776.

We mention other interesting works by Fuss on spherical geometry.
 
  In his memoir 
    \emph{Problematum quorundam sphaericorum solutio} (Solution to certain spherical problems) \cite{Fuss-Problematum-T}, Fuss solves the following three problems:

For a triangle having a fixed base, and whose third vertex is on a given great circle of the sphere, locate this vertex such that:
\begin{enumerate}
\item the angle at that vertex is maximal;
\item the sum of the lengths of the two sides which contain this angle is minimal;
\item the area is maximal.
\end{enumerate}

   In his solution to the first problem, Fuss established a third degree equation from which he obtained conditions for which the problem admits three solutions. In the special case where the angle between the two circles is right, the equation becomes a second degree equation. In the solution of the two other problems, Fuss uses an argument involving an infinitesimal variation of the third vertex. It is interesting to note that all these problems are typical of the field of the calculus of variations, even though Fuss did not use the methods of that theory.

The memoir \emph{De proprietatibus quibusdam ellipseos in superficie sphaerica descriptae} (On some properties of an ellipse traced on a spherical surface) \cite{Fuss-proprietatibus} is in some sense a sequel to the previous one. Fuss develops the theory of ellipses on the sphere. In the same memoir, Fuss studies properties of triangles whose base is fixed and whose vertices are on an ellipse having the vertices of this base as foci. It is possible that this is the problem which motivated Fuss to study ellipses on the sphere. This problem is a variant of the problem studied in the preceding memoir, which consists in finding the properties of triangles with fixed base and whose vertex is on a great circle, with the sum of the lengths of the two sides minimal.

We now talk about Anders Johan Lexell (1740--1784).

Lexell was a Finnish-Swedish astronomer,  mathematician, and physicist. Born in  \AA bo (Finland),\footnote{Until 1812, \AA bo was the capital of Finland, which was part of the Kingdom of Sweden. Today, the Swedish name \AA bo has been replaced by the Finnish name Turku, and Turku is now the sixth largest city in Finland.} he obtained after his studies a position of mathematics professor at the Naval Academy of Uppsala. In 1768, attracted by the presence of Euler in Saint Petersburg, he sent to that academy a memoir titled \emph{Methodus integrandi, nonnullis \ae quationum exemplis illustrata} (Integration methods, illustrated by some examples of equations). We learn from the section on the History of the Academy of Sciences in \cite{Lexell-A} that Euler, who was in charge of examining the manuscript, was highly impressed.  The Academy records also mention that the count Wolodomir Orlov, who was the director of the Academy, after he saw Euler's judgement, objected that the manuscript might well have been written by some clever geometer who was willing to promote Lexell. Euler  replied that among all the working geometers in the world, nobody else than himself or d'Alembert would have been capable of writing such a paper, and that Lexell was unknown to both of them. After this, Orlov sent immediately to Lexell an invitation for an adjunct position  at the academy, which Lexell accepted promptly. 

Lexell was appointed professor of astronomy at the Academy of Sciences of Saint Petersburg in 1771. He became a close collaborator of Euler, with which he started to work mainly on astronomy. He is known for his study of the motion of comets. The king of Sweden offered him in 1775 a position of professor of mathematics at the university of \AA bo, with a permission to stay three more years in Saint Petersburg. Lexell accepted, and the permission was renewed two more times, each time for one year, until 1780. He then returned to his hometown but he stayed there only one year, and he came back to Saint Petersburg in 1781. 

Like Fuss, Lexell was present on the day of Euler's death; the three men were discussing the orbit of Uranus  which had been discovered two years before by William Herschel (1738-1822), and about other scientific matters,  when Euler had his brain attack. Later, Lexell became the first to compute the orbit of Uranus, and his calculations showed that it was a planet rather than a comet, and from perturbations in its orbit, he conjectured the presence of another planet, namely, Neptune, which was discovered subsequently. A comet that has been named after Lexell is famous for having been the closest to Earth in the recorded history.

Lexell was appointed Professor of Mathematics at the Saint Petersburg Academy after Euler's death. Lexell himself got ill and died the year after, at the age of 40.

Lexell worked on a topic which was called ``polygonometry", a generalization of trigonometry where one studies relations between sides and angles of polygons. He obtained formulae for relations between angles and side lengths for quadrilaterals, pentagons, hexagons and heptagons.  He contributed to spherical trigonometry with new and interesting solutions which he took as a basis for his research of comet and planet motion. His name is attached to a theorem on spherical triangles which we discussed in \S \ref{s:Lexell}. Lexell was one of the most prolific members of the Russian Academy of Sciences at his time. He published 66 memoirs in 16 years. We can read in the Proceedings of the Academy of Sciences of Saint Petersburg \cite{Lexell-A} that besides being member of that Academy, he was member of the Royal Academies of Stockholm, of Uppsala and of Turin, and a corresponding member of the Royal Academy of Paris.

 Lexell wrote three papers on spherical trigonometry. We already mentioned his result describing the locus of the vertex of a spherical triangle having a given base and a given area. The result was published in 1784, the year of his death \cite{Lexell-Solutio}. In his two other memoirs, \emph{De proprietatibus circulorum in superficie sphaerica descriptorum} (On the properties of circles traced on the surface of the sphere) \cite{Lexell-proprietatibus} and \emph{Demonstratio nonnullorum theorematum ex doctrina sphaerica} (Proofs of certain theorems according to the spherical doctrine) \cite{Lexell-Demonstratio}, both published in 1782, Lexell obtained several results on spherical geometry, including formulae for the radii of circles inscribed in or circumscribed to spherical triangles and quadrilaterals, and he derived several other results, including a spherical analogue of Heron's formula for the area of Euclidean triangles and a spherical analogue of Ptolemy's theorem.

Among Euler's direct followers who worked on spherical geometry is Friedrich Theodor von Schubert (1758-1825). Like Euler, he was the son of a pastor, and his parents first wanted him to study theology.  He did not follow this path, and he decided to study mathematics and astronomy without teachers and he eventually taught these subjects abroad, as a private teacher. He traveled to Sweden in 1779,  then turned back to Germany, and then moved to Estonia. He was appointed assistant at the Saint Petersburg Academy of Sciences at the class of geography in 1785, that is, two years after the death of Euler, and he became full member of the Academy in 1789. In 1803, he became the director of the astronomical observatory of the Academy. 

Schubert wrote several papers on spherical geometry, most of them related to astronomy or to geography. He was elected member of the academies of Stockholm, Copenhagen, Uppsala, Boston, and others. He is the author of a famous treatise on astronomy in 3 volumes \cite{Schub-Traite}. From 1810 until his death, he was the editor of the German language periodical \emph{Saint Petersburg Zeitung}.

  Schubert also continued the works of Lexell and Fuss. His results include the determination of the loci of the vertices of triangles satisfying some given condition, like the problems solved by Lexell and Fuss. He wrote a paper in which he showed that spherical trigonometry can be developed based on Ptolemy's theorem. Several developments of these works are mentioned in the book of Chasles \cite{Chasles-Apercu}.

In conclusion, spherical geometry is a beautiful subject with many facets and it is related to several other fields of mathematics. Besides presenting part of this subject and the work of Euler and his followers, we hope that this paper can motivate the working mathematician to read the original sources.

 \end{document}